\documentclass[joc]{ipart}
\usepackage{fullpage}
\usepackage{amssymb,amsmath,amsfonts, amscd,url}
\usepackage{array}

\newtheorem{theorem}{Theorem}

\newtheorem{note}[theorem]{Note}
\newtheorem{defn}[theorem]{Definition}
\newtheorem{rem}[theorem]{Remark}

\newtheorem{thm}[theorem]{Theorem}
\newtheorem{cor}[theorem]{Corollary}

\newtheorem{fact}[theorem]{Fact}
\newtheorem{algorithm}[theorem]{Algorithm}
\newcommand{\Forb}{{\rm Forb}}

\newcommand{\wchi}[1]{\chi^{{\rm wk}}_{#1}}
\newcommand{\schi}[1]{\chi^{{\rm st}}_{#1}}
\newcommand{\wchidir}[1]{\chi^{{\rm wk, dir}}_{#1}}
\newcommand{\schidir}[1]{\chi^{{\rm st, dir}}_{#1}}
\newcommand{\chidir}[1]{\chi^{{\rm dir}}_{#1}}

\newcommand{\nonarrow}{\bigcirc}
\newcommand{\unarrow}{-}

\newcommand{\dist}{{\rm dist}}
\newcommand{\forb}{{\rm Forb}}
\newcommand{\arrows}{\mapsto}

\newcommand{\M}{{\bf M}}

\newcommand{\p}{{\bf p}}
\newcommand{\q}{{\bf q}}
\newcommand{\J}{{\bf J}}

\newcommand{\one}{{\bf 1}}
\newcommand{\zero}{{\bf 0}}
\newcommand{\bfd}{{\bf d}}

\newcommand{\hh}{\mathcal{H}}
\newcommand{\F}{\mathcal{F}}
\newcommand{\K}{\mathcal{K}}

\newcommand{\A}{\mathcal{A}}

\newcommand{\chib}{\chi_2}

\newcommand{\E}{{\rm\bf E}}
\newenvironment{proofcite}[1]{\noindent{\bf Proof of #1.\,}}{\hfill$\Box$ \\}

\newcommand{\textdef}{\textit}
\newcommand{\textindef}{\textbf}

\newcommand{\ds}{\displaystyle}

\newcommand{\Aset}{\overleftrightarrow{\mathcal{A}}}
\newcommand{\Pset}{\mathcal{P}}

\begin{document}

\begin{frontmatter}
\title{Multicolor and directed edit distance}
\runtitle{Multicolor and directed edit distance}
\begin{aug}
\author{Maria Axenovich\thanks{This author's research partially supported by NSF grant DMS-0901008 and NSA grant H-98230-09-1-0063.}\ead[label=e1]{axenovic@iastate.edu}}
\address{Department of Mathematics, Iowa State University, Ames, Iowa 50011}
%\email{axenovic@iastate.edu}

\author{Ryan R. Martin\thanks{This author's research partially supported by NSF grant DMS-0901008 and by an Iowa State University Faculty Professional Development grant.}\ead[label=e2]{rymartin@iastate.edu}}
\address{Department of Mathematics, Iowa State University, Ames, Iowa 50011}
%\email{rymartin@iastate.edu}
%\subjclass[2010]{Primary 05C35; Secondary 05C80}
%\keywords{edit distance, hereditary properties, localization, split graphs, colored regularity graphs}
\end{aug}

\maketitle

\begin{abstract}
The editing of a combinatorial object is the alteration of some of its elements such that the
resulting object satisfies a certain fixed property. The edit problem for graphs, when the edges are added or deleted, was first studied independently by the authors and K\'ezdy~\cite{AKM}
%[\textit{J. Graph Theory} (2008) \textbf{58}(2), 123--138]
and by Alon and Stav~\cite{AS1}.
%[\textit{Random Structures Algorithms} (2008) \textbf{33}(1), 87--104]
In this paper, a generalization of graph editing is considered for multicolorings of the complete graph as well as for directed graphs. Specifically, the number of edge-recolorings sufficient to be performed on any edge-colored complete graph to satisfy a given hereditary property is investigated.  The theory for computing the edit distance is extended using random structures and so-called types or colored homomorphisms of graphs.
\end{abstract}
\begin{keyword}[class=AMS]
\kwd[Primary ]{05C35}
\kwd[; Secondary ]{05C80}
\end{keyword}
\begin{keyword}
\kwd{edit distance}
\kwd{hereditary properties}
\kwd{localization}
\kwd{split graphs}
\kwd{colored regularity graphs}
\kwd{\LaTeX}
\end{keyword}
\end{frontmatter}

%\tableofcontents

\section{Introduction}~\\
The combinatorial editing problem is, in general, the problem of finding the smallest number of element-changes such that the resulting combinatorial object satisfies a certain fixed property. The simplest class of objects for which the editing problem was considered is a set of sequences. In fact, the first detailed algorithmic study of editing was motivated by bioinformatics, where sequences over finite alphabets are considered and editing corresponds to changes of the elements in the sequence depicting the mutations in biomolecules. When the desired property consists of a single sequence, studying editing corresponds to investigating the Hamming distance between sequences.   The notion of graph editing was introduced by the authors and K\'ezdy~\cite{AKM} and independently by Alon and Stav~\cite{AS1}. The question considered was: ``How many edges does one need to add or delete in a given graph, such that the result belongs to a given class of graphs?'' The authors showed in~\cite{AKM}, that the answer to this question for hereditary classes could be expressed in terms of the so-called binary chromatic number (also called the colouring number) of the family. Alon and Stav~\cite{AS1} showed that the largest distance from a hereditary property is achieved, asymptotically, by an Erd\H{o}s-R\'enyi random graph.

In this paper, the generalized theory is developed for editing of edge-colored complete graphs and digraphs.  The main result for edge-colored graphs, Theorem~\ref{thm:multicol:easybounds}, is in terms of two parameters: the so-called weak and strong $r$-ary chromatic numbers.  The main result for directed graphs, Theorem~\ref{thm:digraphs:easybounds} is in terms of two parameters: the weak and strong directed chromatic numbers. In each case, the results come from more general theorems, Theorems~\ref{thm:multicol:basics} and~\ref{thm:digraphs:basics} respectively, which deal with generalizing the graph notion of types for the above combinatorial objects.  The analysis is based on using a version of Szemer\'edi's regularity lemma, which we state as Theorem~\ref{thm:multicol:genreglem} (see~\cite{AM} for a proof of Theorem~\ref{thm:multicol:genreglem}), and applying it to an Erd\H{o}s-R\'enyi-type random edge-colored graph or random digraph, respectively.
General bounds on the edit distance function are given, as well as some editing algorithms and computing methods, all of which result from Theorems~\ref{thm:multicol:basics} and~\ref{thm:digraphs:basics}.

The paper is structured as follows. Section~\ref{sec:multicol} deals with the case of multicolorings of the edges of complete graphs.  Section~\ref{sec:digraphs} deals with the case of directed graphs.  In each of these sections we provide definitions, editing algorithms, examples as well as some general theory on the edit distance function.  Most proofs are presented in Subsection~\ref{sec:multicol:proofs} and in Subsection~\ref{sec:digraphs:proofs}.

\section{Multicolorings of the complete graph}~\\
\label{sec:multicol}
\subsection{Basic definitions}~\\
An \textdef{equipartition} of a finite set is a partition in which each pair of partite sets differ in size by at most one.

For a complete graph on vertex set $V$, and a finite set $Q$, we shall say that a \textdef{$Q$-coloring}, or more specifically, a \textdef{$Q$-edge-coloring} of this graph is a pair $G=(V,c)$, where $c: \binom{V}{2} \rightarrow Q$. Since it is sufficient to let $Q=\{1,\ldots,r\}$ for some integer $r$, we will refer to an $\{1,\ldots,r\}$-edge coloring of a complete graph as simply an \textdef{$r$-graph}. For any $r$-graph $G$, disjoint vertex sets $V_i$ and $V_j$ and color $\rho$, $\rho\in\{1,\ldots,r\}$, the expression $E_{\rho}(V_i)$ denotes the set of edges colored $\rho$ with both endpoints in $G[V_i]$ and $E_{\rho}(V_i,V_i)$ denotes the set of edges colored $\rho$ with one endpoint in $V_i$ and the other in $V_j$.  The \textdef{density vector of $V_i$} is an $r$-dimensional vector $\p=(p_1,\ldots,p_r)$, where $p_{\rho}=|E_{\rho}(V_i)|/{\textstyle {|V_i|\choose 2}}$ for $\rho=1,\ldots,r$. The \textdef{density vector of the pair $(V_i,V_j)$} is an $r$-dimensional vector $\p=(p_1,\ldots,p_r)$, where $p_{\rho}=|E_{\rho}(V_i,V_j)|/(|V_i||V_j|)$ for $\rho=1,\ldots,r$. Note that for such density vectors, $\sum_{\rho}p_{\rho}=1$.

In this setting, a \textdef{graph property} is merely a set of $r$-graphs for some positive integer $r\geq 2$. If $G=(V,c)$ and $G'=(V,c')$ are $r$-graphs on $n$ labeled vertices, then
$$\dist(G,G')$$ is the proportion of edges on which the colors differ,
i.e., the number of edges on which the colors in $G$ and $G'$ differ, divided by $\binom{n}{2}$.  We may call this the \textdef{normalized edit distance} between $G$ and $G'$.

For any property $\hh$, a coloring $G$, an integer $n$, we define $\dist(G, \hh)$, $\dist(n, \hh)$, and $\dist(\hh)$ as follows:
\begin{eqnarray*}
   \dist(G, \hh) & := & \min\left\{\dist(G,G') : V(G')=V(G), G'\in \hh\right\} , \\
   \dist(n, \hh) & := & \max\{\dist(G, \hh) : |V(G)|=n\} , \\
   \dist(\hh) & := & \lim_{n\rightarrow\infty}\dist(n,\hh) .
\end{eqnarray*}
Note that $\dist(G,G'),\dist(G,\hh),\dist(n,\hh),\dist(\hh)\in [0,1]$.

The last parameter $\dist(\hh)$ is the limit of the largest proportion of the edges necessary to be changed in a coloring of a complete graph bring it to a property $\hh$; we show the existence of this limit later.

A \textdef{hereditary property of $r$-graphs} (or, simply, \textdef{hereditary property}, where the context is understood) is a set of $r$-graphs that is closed under vertex-deletion and isomorphisms. Let an $r$-graph $G'$ be an \textdef{induced coloring} of an $r$-graph $G$ if $G'$ can be obtained from $G$ by vertex-deletion.

For an $r$-graph, $H$, the family $\forb(H)$ consists of all $r$-graphs that have no (induced) copies of $H$. For every hereditary property, $\hh$, there is a family, $\F(\hh)$, of $r$-graphs such that $\hh=\bigcap_{H\in\F(\hh)}\forb(H)$.  If $\F$ is a family of $r$-graphs, then we use $\forb(\F)$ to denote $\bigcap_{H\in\F}\forb(H)$.~\\

\subsection{The $r$-ary chromatic numbers}~\\
\begin{defn}
For a hereditary property $\hh = \bigcap_{H\in \F(\hh)} \forb(H)$ of $r$-graphs,  a \textindef{weakly-good tuple}  $(a_1, \ldots, a_r)$  is an $r$-tuple of non-negative integers such that for some $H\in \F(\hh)$,  the vertex set $V(H)$ can be partitioned into sets $S_1,\ldots,S_r$ such that, for each $i\in\{1,\ldots,r\}$ with $a_i\neq 0$, the partition can be further refined $S_i=V_{i,1}\cup\cdots\cup V_{i,a_i}$ such that each $V_{i,j}\in S_i$ \textbf{does not} induce an edge of color $i$. The \textindef{weak clique spectrum} of $\hh$ is the set of all tuples $(a_1,\ldots,a_r)$ that are NOT weakly-good. The \textindef{weak $r$-ary chromatic number} of $\hh$, $\wchi{r}(\hh)$, is the maximum $\ell+1$ such that for some non-negative integers $a_1,\ldots,a_r$ with $a_1+\cdots+a_r = \ell$, the tuple $(a_1,\ldots,a_r)$ is in the weak clique spectrum of $\hh$.~\\

For a hereditary property $\hh$,  a \textindef{strongly-good tuple}  $(a_1, \ldots, a_r)$  is an $r$-tuple of non-negative integers such that for some $H\in \F(\hh)$,  the vertex set $V(H)$ can be partitioned into sets $S_1,\ldots,S_r$ such that, for each $i\in\{1,\ldots,r\}$ with $a_i\neq 0$, the partition can be further refined $S_i=V_{i,1}\cup\cdots\cup V_{i,a_i}$ such that each $V_{i,j}\in S_i$ has \textbf{all} edges of color $i$. The \textindef{strong clique spectrum} of $\hh$ is the set of all tuples $(a_1,\ldots,a_r)$ that are NOT strongly-good. The \textindef{strong $r$-ary chromatic number} of $\hh$, $\schi{r}(\hh)$, is the maximum $\ell+1$ such that for some non-negative integers $a_1,\ldots,a_r$ with $a_1+\cdots+a_r = \ell$, the tuple $(a_1,\ldots,a_r)$ is in the weak clique spectrum of $\hh$.

If $\hh=\forb(H)$, then we denote $\wchi{r}(H)=\wchi{r}(\hh)$ and $\schi{r}(H)=\schi{r}(\hh)$.~\\
\end{defn}

\begin{rem}~\\
\begin{itemize}
   \item The weak [strong] clique spectrum is a downset in the partially ordered set of $r$-tuples ordered coordinatewise.  That is, if $(a_1,\ldots,a_r)$ is in the weak [strong] clique spectrum and $(a_1',\ldots,a_r')$ has the property that $0\leq a_i'\leq a_i$ for $i=1,\ldots,r$, then $(a_1',\ldots,a_r')$ is also in that weak [strong] clique spectrum.
   \item Informally, we can partition $V(H)$ into $\wchi{r}(H)$ pieces in which the absent color in each piece is arbitrary, but there is some specification of absent colors for which a $\wchi{r}(H)-1$ piece partition is not possible.
   \item Similarly, we can partition $V(H)$ into $\schi{r}(H)$ pieces in which the required color in each piece is arbitrary, but there is some specification of required colors for which a $\schi{r}(H)-1$ piece partition is not possible.
   \item For any $r\geq 2$ and any hereditary property of $r$-graphs, $\hh$, $\wchi{r}(\hh)\leq\schi{r}(\hh)$.
   \item In the case of $r=2$, notions of strong and weak colorings are identical. Further, if $\hh=\forb(H)$, $\chi_2(H)$ corresponds exactly to the \textdef{binary chromatic number of $H$}, introduced in \cite{AKM}. This is also called the ``colouring number'' in related literature such as Bollob\'as and Thomason~\cite{BT1995,BT2000}.~\\
% ~\footnote{Added this remark and removed this: "Note that if $(a_1,\ldots,a_r)$ is in the clique spectrum of $\hh$ and $G$ is an $r$-graph partitioned into $a_1$ cliques with no edge of color 1, $a_2$ cliques with no edge of color 2, etc., then $G\in\hh$ since $G$ contains no forbidden subgraphs from $\F(\hh)$." and "In the case where $\hh$ is a principal hereditary property, i.e., $\hh=\forb(H)$ for some $H$, then we denote the $r$-ary chromatic number of $\hh$ to be, simply $\chi_r(H)$." and "Note that there is a monotonicity to the clique spectrum.  Precisely: if $(a_1,\ldots,a_r)$ is in a clique spectrum and $a_i'\leq a_i$, $\forall i\in\{1,\ldots,r\}$, then $(a_1',\ldots,a_r')$ is in that clique spectrum.  Similarly, if $(a_1',\ldots,a_r')$ is a good $r$-tuple and $a_i'\leq a_i$, $\forall i\in\{1,\ldots,r\}$, then $(a_1,\ldots,a_r)$ is also a good $r$-tuple. Note also that the zero vector is always in the clique spectrum." -RRM}
\end{itemize}
\end{rem}

\subsubsection{Examples illustrating the $r$-ary chromatic numbers of a hereditary family}~\\

\noindent\textbf{(1)} Let $r=3$ and $\hh$ be a family of $\{1,2,3\}$-colored complete graphs not containing  a triangle $H_1$ with colors $1,1,2$  on its edges and not containing a triangle $H_2$ with colors $2,2,3$ on its edges. So, $\F(\hh)= \{H_1, H_2\}$.

First, the weak $3$-ary chromatic number: Since $r=3$, and $\F(\hh)$ contains a triangle, any  $3$-tuple  $(a_1,a_2,a_3)$ with $a_1+a_2+a_3\geq 3$ must be weakly-good. Indeed, each of $H_1$ and $H_2$ can be vertex-partitioned into three parts such that each part is a single vertex, thus not inducing edges of any colors. Thus, it is sufficient to consider the tuples with $a_1+a_2+a_3 \leq 2$.
%:  $(2,0,0)$,  $(0, 2, 0)$, $(0,0,2)$, $ (1,1, 0)$,  $(1, 0, 1)$, $(0,1,1)$,  $(0, 0, 1)$, $(0, 1, 0)$, $(1, 0, 0)$.
The tuple $(1,0,0)$ is weakly-good since we can partition the vertex set of $H_2$ in one part not containing  edges of color $1$.  Similarly, $(0,0,1)$ is good. By monotonicity, all tuples $(a_1,a_2,a_3) $ with $a_1\geq 1$ or $a_3\geq 1$ are weakly-good.

However $(0,1,0)$ is not weakly-good because both $H_1$ and $H_2$ contain edges of color $2$.  But $(0,2,0)$ is weakly-good  because $H_1$ can be vertex-partitioned in two parts not containing edges of color $2$. Thus, the weak clique spectrum of $\hh$ is $\{(0,1,0),(0,0,0)\}$. For the weak $3$-ary chromatic number, $\wchi{3}(\hh)=2$.

Second, the strong $3$-ary chromatic number. Similar to the above, if $a_1+a_2+a_3\geq 2$, then either $H_1$ or $H_2$ can be partitioned into two parts such that one part is an edge of a specified color and the other part is a vertex.  If $a_1+a_2+a_3\leq 1$, then neither $H_1$ nor $H_2$ can be partitioned into a single monochromatic clique.  Thus, the the strong clique spectrum of $\hh$ is $\{(1,0,0),(0,1,0),(0,0,1),(0,0,0)\}$ and for the strong $3$-ary chromatic number, $\schi{3}(\hh)=2$ also.~\\

\noindent\textbf{(2)} Let $r=3$ and $\hh$ be a family of $\{1,2,3\}$-colored complete graphs not containing a triangle $H_1$ with colors $1,1,2$ on its edges. So, $\F({\hh})=\{H_1\}$. If we follow the previous example, it is easy to see that $(a_1,a_2,1)$ is weekly good for all $a_1,a_2\geq 0$.  Moreover, $(a_1,a_2,0)$ is weakly good as long as $a_1+a_2\geq 2$. Hence, the weak clique spectrum of $\hh$ is $\{(1,0,0),(0,1,0),(0,0,0)\}$ and $\wchi{3}(\hh) = 2$.

For the strong clique spectrum, it is easy to see that $(0,0,2)$ is in that spectrum, but if $a_1+a_2+a_3\geq 3$, then $(a_1,a_2,a_3)$ is strongly-good.  Thus, $\schi{3}(\hh)=3$.~\\

\noindent\textbf{(3)} Let $r=2$, which we can consider to be the graph case. As we have observed, we may disregard the notions of ``weak'' and ``strong'' in our terminology. Let $H$ be a $K_5$ colored with edges colored with colors $1$ and $2$ such that each color class is a $5$-cycle. Let $\hh$ be a family of colorings not containing $H$, i.e., $\F(\hh)=\{H\}$.

We need only consider $2$-tuples (i.e., pairs) $(a_1,a_2)$ with $a_1+a_2\leq 4$. It is relatively easy to see that $(2,1)$, $(1,2)$,  $(3,0)$ and $(0,3)$ are good. The pairs $(2,0)$, $(0,2)$ and $(1,1)$ are not good since $H$ has no monochromatic clique on more than 2 vertices, but has a total of 5 vertices. By monotonicity, $(1,0)$ and $(0,1)$ are also not good. Thus the clique spectrum of $\hh$ is $\{(2,0),(1,0),(1,1),(0,2),(0,1),(0,0)\}$, and $\chi_2(\hh)=3$.~\\

\subsection{A simple editing algorithm}~\\
\label{sec:multicol:simplealg}
Let $\hh$ be a hereditary property of $r$-graphs, such that $\hh = \bigcap_{H\in\F(\hh)}\forb(H)$. Further, let $\ell=\wchi{r}(\hh)-1$ and $(a_1,\ldots,a_r)$ be in the weak clique spectrum with $\sum_{i=1}^ra_i=\ell$.
% and no graph in $\F(\hh)$ can be partitioned into $\ell$ monochromatic cliques with $a_i$ of these cliques of color $i$, for $i=1,\ldots,r$.

Partition $V$ into $r$ sets $S_1,\ldots,S_r$ and further refine the partition such that $S_i=V_{i,1}\cup\cdots\cup V_{i,a_i}$, for $i=1,\ldots,r$ and then recolor the edges in each $V_{i,j}$ by recoloring the edges of color $i$ with some other arbitrary color. This new coloring does not contain any $H\in\F(\hh)$, otherwise the tuple $(a_1,\ldots,a_r)$ would be good for some $H$.

If the sizes of the $V_{i,j}$-s differ by at most one; i.e., $\lfloor n/\ell\rfloor\leq |V_{i,j}|\leq\lceil n/\ell\rceil$, then the number of changes provided by this algorithm is at most $\ell\binom{\lceil n/\ell\rceil}{2}$. Thus,
$$ \dist(\hh) \leq \lim_{n\rightarrow\infty}\frac{\ell\binom{\lceil n/\ell\rceil}{2}}{\textstyle \binom{n}{2}} = \frac{1}{\ell}= \frac{1}{\wchi{r}(\hh)-1}.$$~\\

\subsection{Previous results and new main results}~\\
In \cite{AKM}, the authors and K\'ezdy provide a general bound for $\dist({\hh})$ in the 2-color case.
\begin{theorem}[\cite{AKM}]\label{thm:AKM} For any hereditary property of graphs, $\hh$, with binary chromatic number $\chib\geq 2$,
$$ \frac{1}{2(\chib-1)}\leq \dist({\hh})\leq\frac{1}{\chib-1} . $$
Furthermore, if $\hh=\forb(H)$ such that $H$ is self-complementary, then $\dist(\hh)=\frac{1}{2(\chib(H)-1)}$.
\end{theorem}
\noindent
Here, we show a similar result in the general case.
\begin{theorem}\label{thm:multicol:easybounds}
Let $\hh$ be a hereditary property of $\{1,\ldots,r\}$-edge-colorings of complete graphs.
Let $\wchi{r}=\wchi{r}(\hh)\geq 2$ and $\schi{r}=\schi{r}(\hh)\geq 2$ be the weak and strong (respectively) $r$-ary chromatic numbers of $\hh$. Then,
$$ \frac{1}{r(\schi{r}-1)}\leq\dist(\hh)\leq\frac{1}{\wchi{r}-1} . $$
Furthermore, if $\hh=\forb(H)$ such that all color classes of $H$ are isomorphic, then $\dist(\hh)\leq\frac{1}{r(\wchi{r}(H)-1)}$.
\end{theorem}

We prove Theorem~\ref{thm:multicol:easybounds} in Section~\ref{sec:multicol:proofs}. The upper bound is found in the simple editing algorithm, but to get the lower bound, we need a more general theory.  This is Theorem~\ref{thm:multicol:basics} which is stated in Section~\ref{sec:multicol:algorithm}.  We also prove the result for symmetric colorings in Corollary~\ref{cor:multicol:symm}.  Theorem~\ref{thm:multicol:basics} gives the basic results that deal with computing the edit distance for given hereditary properties. To state these results, we need to provide some preliminary material.~\\

\subsection{The edit distance function}~\\
\subsubsection{Preliminary definitions}
For an $r$-graph, $G=(V,c)$, and some color $\rho\in\{1,\ldots,r\}$, let $E_\rho(G)$ denote the graph on vertex set $V$ corresponding to the edges with color $\rho$ in $c$. For a positive integer $r$, recall that a density vector $\p=(p_1,\ldots,p_r)$ (we also refer to it as a \textdef{probability vector}) is a nonnegative real vector with the property that $\sum_{\rho=1}^rp_i=1$. For any density vector $\p=(p_1,\ldots,p_r)$, and integer $n$, we denote\footnote{Formally, the sizes of the partitions of the edge set should be integral, so we can take the floor function for the sizes of, say $E_{\rho}$ for $\rho=1,\ldots,r-1$ and the size of $E_{r}$ is what remains. Since we fix $p_{\rho}$ for $\rho=1,\ldots,r$ and let $n$ approach infinity, this will make no appreciable difference.}
$$ \dist_n(\p,\hh)=\max\left\{\dist(G,\hh):|V(G)|=n\mbox{ and }|E_\rho(G)|=p_\rho{\textstyle\binom{n}{2}}, \rho=1,\ldots,r \right\} . $$
In Theorem~\ref{thm:multicol:basics}, we show that the following limit exits, which we call the \textdef{edit distance function}:
$$ \dist(\p, {\hh})=\lim_{n\rightarrow\infty}\dist_n(\p,\hh). $$

Having the edit distance function at our disposal, we may also define $\dist(\hh)=\max_{\p}\dist(\p,\hh)$, where the maximum is taken over all density vectors.

\subsubsection{Types of colorings}
\label{sec:multicol:types}

In Section~\ref{sec:multicol:algorithm-analysis}, we define two functions which are described in terms of types of colorings, which allow us to compute the edit distance function. In Section~\ref{sec:multicol:algorithm}, we shall provide an algorithm to do such computation. We define a notion which was called a colored regularity graph (CRG) by Alon and Stav \cite{AS1}, but earlier called a \textdef{type} by Bollob\'as and Thomason \cite{BT1995}. We adopt latter terminology.

\begin{defn}
An \textindef{$r$-type} (or just, \textindef{type}, where the context is understood), $K$, is a pair $(U,\phi)$, where $U$ is a finite set of vertices and $\phi : U\times U \rightarrow 2^{\{1,\ldots,r\}}\setminus \emptyset$, such that $\phi(x,y)=\phi(y,x)$ and $\phi(x,x)\neq\{1,\ldots,r\}$, for all $x, y\in U$. Informally, we will view an $r$-type as a complete graph with a coloring of both vertices and edges using nonempty subsets of $\{1,\ldots,r\}$, where the whole set is a forbidden color on the vertices. The \textindef{sub-$r$-type} of $K$ induced by $W\subseteq U$ is the $r$-type achieved by deleting the vertices $U-W$ from $K$.

We say that an $r$-graph $H=(V,c)$ \textindef{embeds in type $K=(U,\phi)$} if there is a map $\gamma:V\rightarrow U$ such that $c(\{v,v'\})=c_0$ implies $c_0\in\phi(\gamma(v),\gamma(v'))$.
In other words, there is a mapping $\gamma$ that brings each edge of color $c_0$ to a vertex or an edge containing $c_0$ in its color set. If $H$ embeds in type $K$, we write $H\arrows K$, otherwise we write $H \not\arrows K$.
For every hereditary property $\hh$, we let $\K(\hh)$ be the set of all $r$-types such that none of $\F(\hh)$ embeds in that type, i.e.,
$$ \K(\hh)=\left\{K : K\mbox{ is an $r$-type and } H\not\arrows K, \forall H\in\F(\hh)\right\} . $$

% An $r$-graph $G'=(V,c)$ \textindef{has type $K=(\{u_1,\ldots,u_k\},\phi)$} if there is a partition of the vertices of $G'$, $V=V_1\cup\cdots\cup V_k$, such that if $v_i\in V_i$ and $v_j\in V_j$ then $c(v_i,v_j)\in\phi(u_i,u_j)$.

We say that an $r$-graph $G'=(V,c)$ \textindef{has type $K=(U,\phi)$ if} $G'$ embeds into $K$ with mapping $\gamma: V\rightarrow U$ and $\gamma$ is surjective.
\end{defn}

Fact~\ref{fact:multicol:embed} generalizes the ideas underlying the simple editing algorithm in Section~\ref{sec:multicol:simplealg}.
\begin{fact}\label{fact:multicol:embed}
   If $K$ is an $r$-type, $G'$ is of type $K$ and $H$ does not embed into $K$, then $H\not\subseteq G'$.
\end{fact}~\\

\subsection{Editing algorithm using types}~\\
\label{sec:multicol:algorithm}
Let $\p = (p_1, \ldots, p_r)$ and ${\bf w} = (w_1, \ldots, w_k)$ be density vectors; i.e., their entries are nonnegative and sum to 1.  They play different roles, however.  The vector $\p$ will represent a vector of densities, $p_{\rho}$.  That is, the graph $G$ has $p_{\rho}\binom{n}{2}$ edges of color $\rho$. The vector ${\bf w}$ will represent a vector of weights, $w_1,\ldots,w_k$, assigned to the vertices of an $r$-type with vertices $u_1,\ldots,u_k$, respectively.

Let $G=(V,c)$ be an $r$-graph with edges having densities according to the vector $\p=(p_1,\ldots,p_r)$, and $\hh$ be a hereditary property. In order to find an upper bound on $\dist(G, \hh)$, it is sufficient to change $G$ to an $r$-graph, $G'$, such that, for all $H\in\F(\hh)$, $H$ does not embed into the new coloring. In particular, if the resulting coloring has type $K \in \K(\hh)$, then $G'$ is in $\hh$.

\begin{algorithm}\label{algor:multicol:recolor}
Fix a $K=(U, \phi) \in \K(\hh)$ and bring $G$ to a coloring of type $K$ by edge-recoloring. Let $U= \{u_1, \ldots, u_k\}$.
Partition the vertices of $G$ randomly into sets $V_1, \ldots, V_k$ such that the probability of a vertex to be in a part $V_i$ is $w_i$.
Consider an edge $\{x,y\}$ of $G$, let $x\in V_i$, $y\in V_j$, for $i, j \in \{1, \ldots, k\}$. If $c(\{x,y\}) \not\in \phi(\{u_i, u_j\})$, recolor $\{x,y\}$ with a color from $\phi ((u_i, u_j))$. This gives the new $r$-graph $G'$ which, according to Fact~\ref{fact:multicol:embed}, does not admit an embedding of any $H\in\F(\hh)$, thus $G'\in\hh$.
\end{algorithm}

Note that this generalizes the simple algorithm in Section~\ref{sec:multicol:simplealg}. In that algorithm, the type had restricted colorings only on the vertices (possibly of different colors) but each edge receives the color $2^{\{1,\ldots,r\}}$.

\subsubsection{Analysis of the editing algorithm}\label{sec:multicol:algorithm-analysis}
Consider Algorithm~\ref{algor:multicol:recolor} applied with type $K$. Let $G$ be a graph such that the number of edges of color $\rho$ are $p_\rho$ for $\rho=1,\ldots,r$. The expected number of changes is
\begin{eqnarray*}
\E[\#\mbox{ changes}] & = & \binom{n}{2}-\sum_{x,y\in V,~x\neq y} \Pr(\{x,y\}\mbox{ is not changed}) \\
& = & \binom{n}{2}-\sum_{x,y\in V,~x\neq y}\sum _{1\leq i,j \leq k}\Pr(x\in V_i, y\in V_j) {\bf 1}_{c(\{x,y\})\in\phi(u_i,u_j)} \\
& = & \binom{n}{2}-\sum _{1\leq i,j\leq k}w_iw_j \sum_{x,y\in V, ~x\neq y} {\bf 1}_{c(\{x,y\})\in\phi(u_i,u_j)} \\
& = & \binom{n}{2}-\sum _{1\leq i,j\leq k}w_iw_j \sum_{\rho\in\phi(u_i,u_j)} p_\rho\binom{n}{2} \\
\end{eqnarray*}

Let $\M_K(\p)$ be a $k\times k$ matrix such that the $(i,j)$-th entry, $\M_K(\p)(i,j)$, is $1-\sum_{\rho\in\phi(u_i,u_j)}p_\rho$.
Thus, if ${\bf w}=(w_1, \ldots, w_k)$, then
$$ \E[\#\mbox{ changes}]={\bf w}^T\M_K(\p){\bf w}\binom{n}{2} . $$

Finally, we define two functions in terms of the matrix $\M_K(\p)$:
$$ f_K(\p)=\left(\frac{1}{k}\one\right)^T\M_K(\p)\left(\frac{1}{k}\one\right)
\qquad\mbox{and}\qquad
g_K(\p)=\left\{\begin{array}{rrcl}
        \min & \multicolumn{3}{l}{{\bf w}^T\M_K(\p){\bf w}} \\
        \mbox{s.t.} & {\bf w}^T\one & = & 1 \\
        & {\bf w} & \geq & \zero
        \end{array}\right. $$

The $f$ and $g$ functions can be interpreted as follows: If the vertices of an $r$-graph, $G$, are assigned randomly  to parts corresponding to the vertices of $K$, then $f_K(\p)$ and $g_K(\p)$ represent the expectation of the proportion of times that the color of an edge does not map the set of colors in a corresponding vertex or an edge of $K$. The function $f_k(\p)$ is obtained from the uniform distribution, and $g_k(\p)$ is obtained using the optimal distribution $(w_1,\ldots,w_k)$ of the proportion of sizes of parts. Although the $g$ function provides a better bound for $\dist(\p,\hh)$, the linearity of the $f$ function helps prove results from $\dist(\p,\hh)$.~\\

\subsection{Basic results on $r$-graphs}~\\
Theorem~\ref{thm:multicol:basics} summarizes some facts about the edit distance function that generalize easily from results in both and \cite{BM} or \cite{MT}. The proof is in Section~\ref{sec:multicol:basics}. Fix a density vector $\p=(p_1,\ldots,p_r)$.
Formally, the \textdef{random $r$-graph of density $\p$}, or \textdef{random $r$-graph} where the context is clear, is denoted $G(n,\p)$. It is a random variable that is an $\{1,\ldots,r\}$-coloring of the edges of a labeled $K_n$ in which each edge, $e$, is colored independently such that $e$ receives color $\rho$ with probability $p_\rho$.

\begin{theorem}\label{thm:multicol:basics}
Let $\hh$ be a hereditary property of $r$-graphs. Fix an $r$-dimensional density vector $\p$. Then the limit $\dist(\p,{\hh}):=\lim_{n\rightarrow\infty}\dist_n(\p,\hh) $
exists. Moreover,
\begin{enumerate}
\item $\dist(\p,\hh) = \inf_{K\in\K(\hh)}f_K(\p) = \inf_{K\in\K(\hh)}g_K(\p)$; \label{it:multicol:inf}
\item for a fixed $\epsilon>0$, then with probability approaching 1 as $n\rightarrow\infty$,
$$ \dist(\p, \hh)-\epsilon\leq\dist(G(n,\p),\hh) \leq \dist(\p, \hh) ; $$
\item $\dist(\p, \hh)=\lim_{n\rightarrow\infty}{\rm\bf E}[\dist(G(n,\p),\hh)]$;
\item $\dist(\p, \hh)$ is continuous over the domain of $r$-dimensional density vectors and is concave down;\footnote{A function $\psi(\p)$ being concave down means for every pair of density vectors $\p_1,\p_2$ and every real number $t\in[0,1]$, $t\p_1+(1-t)\p_2$ is a density vector and $\psi(t\p_1+(1-t)\p_2)\geq t\psi(\p_1)+(1-t)\psi(\p_2)$.}\label{it:multicol:concon}
\item $\dist(\p,\hh)$ achieves its maximum, $\dist({\hh})$, at some density vector $\p^*_{\hh}$ (in fact, denote the set of all such vectors $\p^*_\hh$) and so,
$$ \dist(\hh)=\lim_{n\rightarrow\infty}{\rm\bf E}[\dist(G(n,\p^*_{\hh}),\hh)] \mbox{; and} $$
\item Both $\p_{\hh}^*$ and $\dist(\hh)$ exist and $\p_{\hh}^*$ is a convex and closed set in $[0,1]^{r-1}$.
\end{enumerate}
\end{theorem}

\begin{rem}
Note that $\p^*_\hh$ typically consists of a single vector, but we abuse notation by denoting the set of such vectors as $\p^*_\hh$ when the vector at which the maximum is obtained is not unique.
\end{rem}

\begin{cor}\label{cor:multicol:symm}
   Let $\hh$ be a symmetric hereditary property; that is, one that has the property such that if the $r$-tuple $(a_1,\ldots,a_r)$ is in the weak clique spectrum of $\hh$, then for any permutation $\varphi$ of $\{1,\ldots,r\}$, the $r$-tuple $(a_{\varphi(1)},\ldots,a_{\varphi(r)})$ is also in the weak clique spectrum.  Then, $$ \dist(\hh)\leq r^{-1}\left(\sum_{i=1}^ra_i\right)^{-1} . $$
   In particular, if $\hh=\forb(H)$ such that all color classes of $H$ are isomorphic, then $\dist(\hh)\leq\frac{1}{r(\wchi{r}(H)-1)}$.
\end{cor}

\begin{proofcite}{Corollary~\ref{cor:multicol:symm}}
   Consider an arbitrary density vector $\p=(p_1,\ldots,p_r)$ and without loss of generality assume that $p_1\leq\cdots\leq p_r$. Choose a permutation of the $a_i$-s such that $a_1\geq\cdots\geq a_r$.  Let $K=(U,\phi)$ be a $r$-type on $\ell=\sum_{i=1}^ra_i$ vertices such that $\phi(u_i,u_j)=\{1,\ldots,r\}$ if $i\neq j$ and there are exactly $a_j$ vertices $u$ such that $\phi(u,u)=\{1,\ldots,r\}-\{j\}$.

   The off-diagonal entries of $\M_K(\p)$ are zero and so it is easy to see that $f_K(p)=\ell^{-2}\sum_{i=1}^ra_ip_i$.  We can use a correlation inequality such as FKG~\cite{FKG} to see that
   $$ f_K(p)=\ell^{-2}\sum_{i=1}^ra_ip_i\leq\ell^{-2}r^{-1}\left(\sum_{i=1}^ra_i\right)\left(\sum_{i=1}^rp_i\right)=r^{-1}\ell^{-1} . $$
   To finish the proof observe that, in the case of $\hh=\forb(H)$, $\ell=\sum_{i=1}^ra_i=\wchi{r}(H)-1$.
\end{proofcite}

\subsection{Example: triangles}~\\
Theorem~\ref{thm:multicol:triangles} gives some basic results on examples of hereditary properties of $r$-graphs defined by triangles.  The proof is in Section~\ref{sec:multicol:triangles}.
\begin{thm}\label{thm:multicol:triangles}
Let $r=3$ and consider hereditary properties of $r$-graphs.
\begin{enumerate}
\item If $\F$ is a family of that consists of a single monochromatic triangle, then $\dist(\Forb(\F))=1/2$.\label{it:multicol:triangles:1}
\item If $\F$ is a family that consists of a single triangle with two edges colored $1$ and one edge colored $2$, then $\dist(\Forb(\F))=1/2$.\label{it:multicol:triangles:2}
\item If $\F$ is a family that consists of two monochromatic triangles of different colors, then $\dist(\Forb(\F))=1/2$.\label{it:multicol:triangles:3}
\item If $\F$ is a family that consists of all six bi-chromatic triangles, then $\dist(\Forb(\F))=2/3$.\label{it:multicol:triangles:4}
\item If $\F$ is a family that consists of a single rainbow triangle, then $\dist(\Forb(\F))=1/3$.\label{it:multicol:triangles:5}
\end{enumerate}
\end{thm}~\\

\subsection{Proofs}~\\
\label{sec:multicol:proofs}
\subsubsection{Proof of Theorem \ref{thm:multicol:easybounds}}
The upper bound for this theorem is proven by the simple editing algorithm from Section~\ref{sec:multicol:simplealg}.

For the lower bound, we apply part (\ref{it:multicol:inf}) of Theorem~\ref{thm:multicol:basics}, which states that $\dist(\p,\hh)=\inf_{K\in\K(\hh)}f_K(\p)$. Consider an arbitrary $K=(U,\phi)\in\K(\hh)$, an $r$-type on $k$ vertices. Let $\tilde{K}$ be an auxiliary graph with vertex set $U$ such that $u$ and $u'$ are adjacent in $\tilde{K}$ if and only if $\phi(u,u')=\{1,\ldots,r\}$. We observe that $\tilde{K}$ has no clique on $\schi{r}=\schi{r}(\hh)$ vertices, otherwise for some $H\in\F(\hh)$, $H\mapsto K$. Using Tur\'an's theorem, the number of edges of $\tilde{K}$ is at most $\frac{\schi{r}-2}{\schi{r}-1}\cdot\frac{k^2}{2}$.

Let $\p=\frac{1}{r}{\bf 1}$. Consider the matrix ${\bf M}_K(\p)$ and observe that every entry is either zero or is a positive integer multiple of $1/r$. The zero entries correspond exactly to pairs with $\phi$ value equal to $\{1,\ldots,r\}$. Thus, this matrix ${\bf M}_K(\p)$ has at least $k^2-2\left(\frac{\schi{r}-2}{\schi{r}-1}\cdot\frac{k^2}{2}\right)\geq\frac{k^2}{\schi{r}-1}$ entries with value at least $1/r$.
Therefore, $f_K(\p)=\frac{1}{k^2}{\bf 1}^T{\bf M}_K(\p){\bf 1}$ is at least $\frac{1}{r(\schi{r}-1)}$. Since $K$ was arbitrary, this gives a lower bound for $\dist(\p,\hh)$.~\\

\subsubsection{Proof of Theorem \ref{thm:multicol:basics}}~\\
\label{sec:multicol:basics}
Let $f(\p)=\inf_{K\in K(\hh)}f_K(\p)$ and let $g(\p)=\inf_{K\in K(\hh)}g_K(\p)$.~\\

%=====================================================
%=====================================================
%=====================================================
\noindent\textbf{A:} Upper bound on $\dist(\p,\hh)$.~\\
%=====================================================
%=====================================================
%=====================================================
Let $G$ be an $r$-graph with the density of its $i$-th color class be $p_{\rho}$ for $\rho=1,\ldots,r$.
Let $K \in \K(\hh)$. Apply the editing algorithm in Section~\ref{sec:multicol:algorithm} to $G$ using $K$. The analysis of the algorithm in Section~\ref{sec:multicol:algorithm-analysis} gives that the expected number of changes is
$f_K(\p)\binom{n}{2}$ and so $\dist_n(\p,\hh)\leq f(\p)\binom{n}{2}$.~\\

%=====================================================
%=====================================================
%=====================================================
\noindent\textbf{B:} Equality of $f$ and $g$.~\\
%=====================================================
%=====================================================
%=====================================================
By the definition of $g_K(\p)$, it is easy to see that $g_K(\p)\leq f_K(\p)$ for every density vector $\p$. Therefore, $g(\p)\leq f(\p)$.  For the other direction, we will use $K$ and its optimal weight vector ${\bf w}^*=\{w_1,\ldots,w_k\}$, where $w_i$ corresponds to $v_i\in V(K)$ in order to construct a sequence of CRGs, $\{K_{\ell}\}$ such that $\lim_{\ell\rightarrow\infty}f_{K_{\ell}}(\p)=g_K(\p)$.

First, choose $\ell$ large enough to ensure that $w_i\ell\geq 2$ for $i=1,\ldots,k$.  Then, for each vertex $u_i\in V(K)$, create $\lfloor w_i\ell\rfloor$ copies of $u_i$ in the following sense:  Let $u_i'$ and $u_j''$ be copies of $u_i$ and $u_j$, respectively, where $u_i,u_j\in V(K)$.  Let $\phi$ be the coloring function of $K$ and $\phi'$ be the coloring function of $K_{\ell}$.  If $i\neq j$, then $\phi'(u_i',u_j'')=\phi(v_i,v_j)$.  If $i=j$ and $v_i'\neq v_i''$, then $\phi'(u_i',u_i'')=\phi(v_i,v_i)$.  Finally, $\phi'(v_i',v_i')=\phi(v_i,v_i)$.

The $(i,j)$-th block is a $\lfloor w_i\ell\rfloor\times\lfloor w_j\ell\rfloor$ matrix and each entry of the $(i,j)$-th block is the same as the $(i,j)$-th entry of ${\bf M}_K(\p)$.

If we denote the $(i,j)$-th entry of ${\bf M}_K(\p)$ by $m_{ij}$, then
$$ \begin{array}{rclcl}\ds
   f_{K_{\ell}}(\p) & = & \ds\frac{1}{|V(K)|^2}\one^T{\bf M}_{K_{\ell}}(\p)\one & = & \ds\left(\sum_i\lfloor w_i\ell\rfloor\right)^{-2}\sum_{i,j}m_{ij}\lfloor w_i\ell\rfloor\lfloor w_j\ell\rfloor \\~\\
   & \leq & \ds\ell^2\left(\sum_i\lfloor w_i\ell\rfloor\right)^{-2}\sum_{i,j}m_{ij}w_iw_j
%   & = & \ds\ell^2\left(\sum_i\lfloor w_i\ell\rfloor\right)^{-2}({\bf w}^*)^T{\bf M}_K(\p){\bf w} \\
   & = & \ds\ell^2\left(\sum_i\lfloor w_i\ell\rfloor\right)^{-2}g_K(\p) \\~\\
   & \leq & \ds\ell^2\left(\sum_i(w_i\ell-1)\right)^{-2}g_K(\p)
   & = & \ds\frac{\ell^2}{(\ell-k)^2}g_K(\p) .
\end{array} $$
Taking $\ell\rightarrow\infty$, we see that $\lim_{\ell\rightarrow\infty}f_{K_{\ell}}(\p)\leq g_K(\p)$. Consequently, for any $K\in\K(\hh)$,
$$ f(\p)=\inf_{\tilde{K}\in\K(\hh)}f_{\tilde{K}}(\p)\leq\lim_{\ell\rightarrow\infty}f_{K_{\ell}}(\p)\leq g_K(\p) . $$
Take the infimum over all $K\in\K(\hh)$, and we have that $f(\p)\leq g(\p)$.~\\

%=====================================================
%=====================================================
%=====================================================
\noindent\textbf{C:} Lower bound on $\dist(\p,\hh)$ using the random graph.~\\
%=====================================================
%=====================================================
%=====================================================
% Here we use the usual relaxation of the different models of random graphs. That is, the distance of the random density $\p$, $r$-graph on $n$ vertices from $\hh$ can be approximated by $G(n,\p)$.
%This approach mirrors that of Alon and Stav~\cite{AS1}.
We apply Theorem~\ref{thm:multicol:reglemapp}, which is given in~\cite{AM}.  Theorem~\ref{thm:multicol:reglemapp} is a corollary of Theorem~\ref{thm:multicol:genreglem}, a relatively straightforward generalization to $r$-graphs and digraphs of a theorem by Alon, Fischer, Krivelevich and M. Szegedy~\cite{AFKSz}, which is suitable for induced graphs.

In an $r$-graph, the \textdef{density vector of a pair of disjoint sets of vertices $(V_i,V_j)$} is simply $\bfd(V_i,V_j):=\left(d_1(V_i,V_j),\ldots,d_r(V_i,V_j)\right)$. So we can state the general version of the regularity lemma. For all definitions of regularity, see~\cite{AFKSz}.
\begin{thm}[Alon, et al.~\cite{AFKSz}]\label{thm:multicol:genreglem}
\newcommand{\Smain}{S_{\ref{thm:multicol:genreglem}}}
\newcommand{\delmain}{\delta_{\ref{thm:multicol:genreglem}}}
   Fix $r\geq 2$.  For every $m$ and function $\mathcal{E}$ with $\mathcal{E}:{\mathbb{N}}\rightarrow (0,1)$, there exist $S=\Smain(r,m,\mathcal{E})$ and $\delta=\delmain(r,m,\mathcal{E})$ with the following property: \\
   \indent If $G$ is a graph [$r$-graph, digraph] with $n\geq S$ vertices then there exist an equipartition $\A=\{V_i : 1\leq i\leq k\}$ of $G$ and an induced subgraph [induced $r$-subgraph, induced subdigraph] $G'$ of $G$, with an equipartition $\A'=\{V_i' : 1\leq i\leq k\}$ of the vertices of $G'$ that satisfy:
   \begin{itemize}
      \item $S\geq k\geq m$.
      \item $V_i'\subset V_i$ for all $i\geq 1$, and $|V_i'|\geq\delta n$.
      \item In the equipartition $\A'$, \textbf{all} pairs are $\mathcal{E}(k)$-regular.
      \item All but at most $\mathcal{E}(0)\binom{k}{2}$ of the pairs $1\leq i<i'\leq k$ are such that $\|\bfd(V_i,V_{i'})-\bfd(V_i',V_{i'}')\|_{\infty}<\mathcal{E}(0)$.
   \end{itemize}
\end{thm}

We use Theorem~\ref{thm:multicol:genreglem} in order to prove Theorem~\ref{thm:multicol:reglemapp}, which is the result that we need.
\begin{thm}[\cite{AM}]\label{thm:multicol:reglemapp}
   Let $G'$ be an $r$-graph in hereditary property $\hh=\bigcap_{H\in\F(\hh)}\forb(H)$ and ${\bf p}=(p_1,\ldots,p_r)$ be a density vector.  Then, there exists an $r$-type $K\in\K(\hh)$ such that $H\not\mapsto K$ for all $H\in\F(\hh)$ and with probability going to $1$ as $n\rightarrow\infty$, $\dist(G_{n,{\bf p}},\hh)\geq f_K({\bf p})\binom{n}{2}-o(n^2)$.
\end{thm}

The proof of Theorem~\ref{thm:multicol:reglemapp} from Theorem~\ref{thm:multicol:genreglem} is straightforward and the details are given in~\cite{AM}.  We begin with $G$ distributed according to $G(n,\p)$ and typical in the sense that any Szemer\'edi partition will have every pair $n^{-0.4}$-regular. Let $G'$ be the graph of smallest distance from $G$ and apply Theorem~\ref{thm:multicol:genreglem}. The resulting partition $\A'$ describes a type $K$ which must be in $\K(\hh)$.  Furthermore, the number of changes required to ensure that $G'$ has partition $\A$ is very close to $f_K(\p)$ because almost every pair in $\A$ has the same density as in $\A'$.

Using part A, we see that for any $\epsilon>0$, with probability approaching 1 as $n\rightarrow\infty$,
\begin{equation}\label{eq:randomlimit}
  f(\p)-\epsilon/2\leq\dist(G(n,\p),\hh)\leq\dist(\p,\hh)\leq f(\p) .
\end{equation}~\\

We can now combine A, B and C. Take the limit of (\ref{eq:randomlimit}) as $n\rightarrow\infty$, and we obtain that for all $\epsilon>0$, $f(\p)-\epsilon/2\leq \dist(\p,\hh)\leq f(\p)$. Hence, $\dist(\p, \hh)=f(\p)=g(\p)$. Moreover, we can replace the second term with ${\rm\bf E}[\dist(G(n,\p),\hh)]$ because that random variable is bounded (in $[0,1]$) and so (\ref{eq:randomlimit}) occurring with high probability implies that the random variable is concentrated around its mean, which approaches $\dist(\p,\hh)$. This verifies parts (1), (2) and (3) of the theorem.~\\

%=====================================================
%=====================================================
%=====================================================
%=====================================================
\noindent\textbf{D:} Continuity of $f$.~\\
%=====================================================
%=====================================================
%=====================================================
Because the set of $r$-types is countable, we can linearly order $\K(\hh)$ to be $K_1,K_2,\ldots$. For every density vector $\p$, set $m_{\ell}(\p)=\min_{i\leq\ell}f_{K_i}(\p)$.

We want to show that each function $m_{\ell}$ is Lipschitz with coefficient 1 with respect to the $L^1$ metric. Let $\p=(p_1,\ldots,p_{\rho})$ and $\q=(q_1,\ldots,q_{\rho})$ be density vectors and define $r$-types $K_{\p},K_{\q}\in\{K_1,\ldots,K_{\ell}\}$ on $k_{\p},k_{\q}$ vertices, respectively, such that $m_{\ell}(\p)=f_{K_{\p}}(\p)$ and $m_{\ell}(\q)=f_{K_{\q}}(\q)$. Then, using the matrix definition of $f$ and the definition of $m_{\ell}$ as a minimum of linear functions,
\begin{eqnarray*}
  f_{K_{\p}}(\p)-f_{K_{\p}}(\q) \leq & f_{K_{\p}}(\p)-f_{K_{\q}}(\q) & \leq f_{K_{\q}}(\p)-f_{K_{\q}}(\q) \\
  \left(\frac{1}{k_{\p}}\one\right)^T\M_{K_{\p}}(\p-\q)\left(\frac{1}{k_{\p}}\one\right) \leq & f_{K_{\p}}(\p)-f_{K_{\q}}(\q) & \leq \left(\frac{1}{k_{\q}}\one\right)^T\M_{K_{\q}}(\p-\q)\left(\frac{1}{k_{\q}}\one\right)
\end{eqnarray*}

Since each of the entries in matrices $\M_{K_{\p}}$ and $\M_{K_{\q}}$ is between zero and one,
and the number of entries in these matrices is $k_{\p}^2$ and $k_{\q}^2$, respectively, it is the case that
$$ \left|f_{K_{\p}}(\p)-f_{K_{\q}}(\q)\right| \leq \|\p-\q\|_1 . $$

Since $\{m_{\ell}\}_{\ell\geq 1}$ is Lipschitz, Definition 7.22 from Rudin~\cite{Rudin} says that the sequence of functions is equicontinuous. The sequence is also pointwise bounded above by $f_{K_1}(\p)$ and below by $0$. By Theorem 7.25(b) from~\cite{Rudin} the sequence $\{m_{\ell}\}_{\ell\geq 1}$ has a uniformly convergent subsequence. Since $\{m_{\ell}\}_{\ell\geq 1}$ is an equicontinuous, each member is itself continuous. Theorem 7.12 from~\cite{Rudin} gives that the aforementioned uniformly convergent subsequence has a continuous limit. The monotonicity of $\{m_{\ell}\}_{\ell\geq 1}$ gives that the limit of any subsequence is the same as the pointwise limit of the sequence itself, namely $\lim_{\ell\rightarrow\infty}m_{\ell}=\inf_{K\in\K(\hh)}f_K=\dist(\hh)$.~\\

%=====================================================
%=====================================================
%=====================================================
\noindent\textbf{E:} Concavity.~\\
%=====================================================
%=====================================================
%=====================================================
Let $\p_1$ and $\p_2$ be density vectors and $t\in[0,1]$ be a real number.  Observe that $t\p_1+(1-t)\p_2$ is still a density vector and, hence, in the domain.  Furthermore,
\begin{eqnarray*}
  f(t\p_1+(1-t)\p_2) & = & \inf_{K\in\K(\hh)}\left\{f_K\left(t\p_1+(1-t)\p_2\right)\right\} \\
  & = & \inf_{K\in\K(\hh)}\left\{tf_K(\p_1)+(1-t)f_K(\p_2)\right\} \\
  & \geq & t\left(\inf_{K\in\K(\hh)}\left\{f_K(\p_1)\right\}\right) +(1-t)\left(\inf_{K\in\K(\hh)}\left\{f_K(\p_2)\right\}\right) \\
  & = & t f(\p_1)+(1-t)f(\p_2) .
\end{eqnarray*}
This gives concavity.~\\

Using D and E, we obtain part (4) directly and the fact that $g_{\hh}$ achieves its maximum follows from continuity (and compactness) and Theorem 4.16 from~\cite{Rudin}. Let $S$ be the set of density vectors $\p$ such that $\dist(\p,\hh)=\dist(\hh)$. The set $S$ must be convex set, because if $\dist(\p_1,\hh)=\dist(\p_2,\hh)=\dist(\hh)$, then by continuity and concavity, the line segment that connects $\p_1$ and $\p_2$ must consist of vectors in $S$. The set $S$ must be closed because a corollary to Theorem 4.8 from~\cite{Rudin} says that, under a continuous mapping, the inverse image of a closed set is closed. Since $\dist(\p,\hh)$ is a continuous function and $S$ is the inverse image of the closed set, $\{\dist(\hh)\}$, then $S$ is closed. This verifies parts (4), (5) and (6) of the theorem and concludes the proof.\hfill$\Box$ \\

\subsubsection{Proof of Theorem~\ref{thm:multicol:triangles}}~\\
\label{sec:multicol:triangles}

\noindent\textbf{(\ref{it:multicol:triangles:1})} In order to destroy all copies of a monochromatic $1$-colored triangle in an arbitrary coloring of $K_n$, it is sufficient to split the vertex set into two parts and recolor all edges within these parts in color $2$. This requires at most $\frac{1}{2}\binom{n}{2}$ changes. To see the lower bound, consider $K_n$ with all edges colored $1$.  After all editing is done to ensure that color class $1$ has no triangles, color class $1$ is triangle-free, having at most $\frac{n^2}{4}$ edges. Thus, at least $\frac{n^2}{4}=\frac{1}{2}\binom{n}{2}+o(n^2)$ edges must have been changed.~\\

\noindent\textbf{(\ref{it:multicol:triangles:2})} In order to destroy all such triangles, it suffices to equipartition the vertex set into two parts and recolor all edges within these parts to color $3$. This requires at most $\frac{1}{2}\binom{n}{2}$ changes. To see the lower bound, consider $K_n$ on vertex set with equipartition $V_1\cup V_2$. Let all edges between $V_1$ and $V_2$ be colored $1$ and let all edges within parts $V_i$, $i=1, 2$ be colored $2$.  We may assume that the only editing operations are recoloring an edge of color $1$ into color $3$ and recoloring an edge of color $2$ into color $3$ because this editing will never create a forbidden triangle.  Let $c$ be such a recoloring not containing triangles with two edges of color $1$ and one edge of color $2$. Let $G$ be an auxiliary graph corresponding to edges of color $3$ in this coloring. The complement of $G$ can not have any triangles with vertices in both $V_1$ and $V_2$. It is easy to prove by induction on $n$ that a graph with satisfying such a condition could have at most $\frac{1}{2}\binom{n}{2}$ edges. Therefore $G$ has at least $\frac{1}{2}\binom{n}{2}$ edges, and this corresponds to the number of changes made.~\\

\noindent\textbf{(\ref{it:multicol:triangles:3})} Assume that $\F$ consists of a triangle with all edges colored $1$ and of a triangle with all edges colored $2$. In order to destroy both of these triangles in an any coloring, as in the previous case, it is sufficient to equipartition the vertex set into two parts and recolor all edges within these parts in color $3$. This requires at most $\frac{1}{2}\binom{n}{2}$ changes.

As to the lower bound, fix $\p=(1/2,1/2,0)$ and consider a $3$-type, $K\in\K(\hh)$, on $k$ vertices. Each of the vertices must have color $3$.  By Tur\'an's theorem, at least $\binom{k}{2}-\lfloor k^2/4\rfloor=\lceil (k^2-2k)/4\rceil$ edges cannot have color $1$ and at least $\lceil (k^2-2k)/4\rceil$ edges cannot have color $2$. Hence, if we consider the off-diagonal entries of $\M_K(\p)$, the sum is at least $\frac{1}{2}\lceil (k^2-2k)/4\rceil+\frac{1}{2}\lceil (k^2-2k)/4\rceil$. So, for any such $K$,
$$ f_K(p)\geq\frac{1}{k^2}\left[k+2\left\lceil\frac{k^2-2k}{4}\right\rceil\right]\geq\frac{1}{2} . $$
As a result, $\inf_{K\in\K(\hh)}f_K(\p)\geq 1/2$.~\\

\noindent\textbf{(\ref{it:multicol:triangles:4})} It is suffices to recolor edges of colors $1$ or $2$ into color $3$. As a result, all forbidden colored triangles will be destroyed via at most $\frac{2}{3}\binom{n}{2}$ changes. In fact, for fixed $\p=(p_1,p_2,p_3)$, at most $(1-\max\{p_1,p_2,p_3\})\binom{n}{2}$ changes suffice.

To see the lower bound, consider a $3$-type $K\in\K(\hh)$ on $k$ vertices. The vertices must be monochromatic and, in addition, the edges incident to a vertex must share the color of that vertex. Otherwise, there would be a bichromatic triangle $H$ with $H\mapsto K$. This implies, however, that $K$ must be entirely monochromatic.  Hence, $g_K(\p)\geq 1-\max\{p_1,p_2,p_3\}$.

Note the this determines not only $\dist(\hh)$, but the entire function $\dist(\p,\hh)=1-\max\{p_1,p_2,p_3\}$.~\\

\noindent\textbf{(\ref{it:multicol:triangles:5})} Observe that in order to destroy all rainbow triangles using colors $1$, $2$ and $3$, it is sufficient to edit the smallest of these color classes, thus performing at most a $\min\{p_1,p_2,p_3\}$ proportion of changes.

For the lower bound, simply observe that no edge in any $K\in\K(\hh)$ can be trichromatic.  Otherwise, that edge, together with any vertex to which it is incident admits a mapping of a rainbow triangle. Hence, each entry of $M_K(p)$ is at least $\min\{p_1,p_2,p_3\}$ and so $f_K(p)\geq\min\{p_1,p_2,p_3\}$.  Hence $\dist(\p,\hh)=\min\{p_1,p_2,p_3\}$ and $\dist(\hh)=1/3$.\hfill$\Box$ \\

%=================================================
%=================================================
%=================================================
%=================================================
%=================================================
%=================================================

\section{Directed graphs}~\\
\label{sec:digraphs}
\subsection{Basic definitions}~\\
We give a number of definitions that are similar to the case of $r$-graphs, however, there are some important distinctions.

\begin{defn}
A \textindef{simple directed graph} or \textindef{digraph} is defined to be a pair $(V,E)$ where $V$ is a labeled vertex set, $E\subseteq (V)_2$ and $(V)_2$ denotes the set $V\times V-\{(v,v) : v\in V\}$. We will also view this as a coloring; that is, a digraph is a pair $(V,c)$ where $c:(V)_2\rightarrow\{\nonarrow,\unarrow,\leftarrow,\rightarrow\}$ is a function which has the property that, for distinct $v,w$,
\begin{itemize}
   \item $c(v,w)=c(w,v)$ if and only if $c(v,w)\in\{\nonarrow,\unarrow\}$ and
   \item $c(v,w)=\rightarrow$ if and only if $c(w,v)=\leftarrow$.
\end{itemize}
Let $\Aset:=\{\nonarrow,\unarrow,\leftarrow,\rightarrow\}$.  Here we interpret the color $c(v,w)=\nonarrow$ to mean that neither $(v,w)$ nor $(w,v)$ are in $E$, the color $c(v,w)=\unarrow$ to mean that both $(v,w)$ and $(w,v)$ are in $E$ and the color $c(v,w)=\rightarrow$ to mean that $(v,w)\in E$ and $(w,v)\not\in E$. \label{defn:digraphs:basic}
\end{defn}

For any digraph $G$ on fixed vertex set $\{v_1,\ldots,v_n\}$, disjoint vertex sets $V_i$ and $V_j$ and color $\rho$, $\rho\in\Aset$, the expression $E_{\rho}(V_i)$ denotes the set of pairs $\{v_i,v_i'\}$ with $v_i,v_i'\in V_i$, $i<i'$ and $c(v_i,v_i')=\rho$. The expression $E_{\rho}(V_i,V_j)$ denotes the set of pairs $\{v_i,v_j\}$ with $v_i\in V_i$, $v_j\in V_J$ and $c(v_i,v_j)=\rho$.  Hence, $E_{\leftarrow}(V_i,V_j)=E_{\rightarrow}(V_j,V_i)$.  As it happens, we will be able to assume, as in the proof of Theorem~\ref{thm:multicol:basics}, that our graphs are random. We will also be able to assume that, among the pairs that have directed edges, a $\leftarrow$ is as likely as $\rightarrow$.  Hence, we will postpone the definition of a density vector for directed graphs.

\begin{defn}
We say that $\Pset\subseteq\Aset$ is a \textindef{palette} if either none or both of ``$\rightarrow$'' and ``$\leftarrow$'' are in $\Pset$. There are 5 possible nontrivial palettes:
\begin{enumerate}
\setcounter{enumi}{-1}
\item $\Pset_0=\Aset$ is the most general case. \label{it:palette:A}
\item $\Pset_{\rm compl}=\{\unarrow,\leftarrow,\rightarrow\}$ is the case of simple digraphs such that every pair of vertices has at least one arc between them. \label{it:palette:noempty}
\item $\Pset_{\rm orien}=\{\nonarrow,\leftarrow,\rightarrow\}$ is the case of \textdef{oriented graphs}; that is, no pair of vertices has two arcs between them. \label{it:palette:nounarrow}
\item $\Pset_{\rm undir}=\{\nonarrow,\unarrow\}$ is the case of simple, undirected graphs.  \label{it:palette:undir}
\item $\Pset_{\rm tourn}=\{\leftarrow,\rightarrow\}$ is the case of \textdef{tournaments}. \label{it:palette:tourn}
\end{enumerate}
\end{defn}

The palette is the universe in which the editing takes place.  That is, if $\nonarrow$ is not in the palette, then no pair $(v,w)$ can be changed to color $\nonarrow$ in the editing process.

If $\Pset$ is a fixed palette and $G=(V,c)$ and $G'=(V,c')$ are digraphs with colors in $\Pset$, then $\dist(G,G')$ is the proportion of edges on which the colors differ; i.e., the number of edges on which the colors differ, divided by $\binom{n}{2}$.\footnote{Here, we can talk about pairs because the color of the pair $(v,w)$ determines the color of the pair $(w,v)$.} We may call this the \textdef{normalized edit distance} between $G$ and $G'$. For any property $\hh$, a simple digraph $G$ with all edge-colors in palette $\Pset$, an integer $n$, we define $\dist(G, \hh)$, $\dist(n, \hh)$, and $\dist(\hh)$ similarly to the  multicolor case.

A \textdef{hereditary property of digraphs with respect to palette $\Pset$} (or, simply, \textdef{hereditary property}, where the context is understood) is a set of digraphs with all edge-colors in $\Pset$ that is closed under vertex-deletion and isomorphisms. Let a digraph $G'$ be an \textdef{induced digraph} of $G$ if $G'$ can be obtained from $G$ by vertex-deletion.
For a fixed palette, $\Pset$ and a digraph, $H$, the family $\forb(H)$ (the palette will be understood) consists of all digraphs with edge-colors in $\Pset$ that have no (induced) copies of $H$. For every palette $\Pset$ and every hereditary property $\hh$ with respect to $\Pset$, there is a family, $\F(\hh)$, of digraphs such that $\hh=\bigcap_{H\in\F(\hh)}\forb(H)$.  If $\F$ is a family of digraphs, then we use $\forb(\F)$ to denote $\bigcap_{H\in\F}\forb(H)$.~\\

\subsection{The directed chromatic numbers}~\\
\begin{defn}
For a hereditary property $\hh=\bigcap_{H\in\F(\hh)}\forb(H)$ and a palette $\Pset$,  a \textindef{weakly-good triple} $(a_0,a_1,a_2)$ is a triple of non-negative integers such that for some $H\in \F(\hh)$, the vertex set $V(H)$ can be partitioned into sets $S_0,S_1,S_2$ such that, for each $i\in\{0,1,2\}$ with $a_i\neq 0$, the partition can be further refined $S_i=V_{i,1}\cup\cdots\cup V_{i,a_i}$ and
\begin{itemize}
   \item each $V_{0,j}$ does not induce a nonedge (i.e., does not induce an edge of color $\nonarrow$),
   \item each $V_{1,j}$ ensures that the directed edges induced by $V_{1,j}$ form an acyclic digraph, and
   \item each $V_{2,j}$ does not induce a bidirectional edge (i.e., does not induce an edge of color $\unarrow$).
\end{itemize}
The \textindef{weak clique spectrum} of $\hh$ with respect to a palette $\Pset$ is the set of all triples $(a_0,a_1,a_2)$ that are NOT weakly-good and such that
$a_0=0$ if $\nonarrow\not\in\Pset$,
$a_1=0$ if $\{\rightarrow,\leftarrow\}\cap\Pset=\emptyset$,
$a_2=0$ if $\unarrow\not\in\Pset$. The \textindef{weak directed chromatic number}, $\wchidir{\Pset}(\hh)$,  of $\hh$ with respect to a palette $\Pset$ is the maximum $\ell +1$ such that for some non-negative integers $a_0,a_1,a_2$, with $a_0+a_1+a_2=\ell$, the triple $(a_0,a_1,a_2)$ is in the weak clique spectrum of $\hh$. We merely use $\wchidir{}(\hh)$ for the weak directed chromatic number if the palette is understood.~\\

For a hereditary property $\hh=\bigcap_{H\in\F(\hh)}\forb(H)$ and a palette $\Pset$,  a \textindef{strongly-good triple} $(a_0,a_1,a_2)$ is a triple of non-negative integers such that for some $H\in \F(\hh)$, the vertex set $V(H)$ can be partitioned into sets $S_0,S_1,S_2$ such that, for each $i\in\{0,1,2\}$ with $a_i\neq 0$, the partition can be further refined $S_i=V_{i,1}\cup\cdots\cup V_{i,a_i}$ and
\begin{itemize}
   \item each $V_{0,j}$ induces only nonedges (i.e., all edges are of color $\nonarrow$),
   \item each $V_{1,j}$ ensures that the directed edges induced by $V_{1,j}$ induce a transitive tournament, and
   \item each $V_{2,j}$ induces only bidirectional edges (i.e., all edges are of color $\unarrow$).
\end{itemize}
The \textindef{strong clique spectrum} of $\hh$ with respect to a palette $\Pset$ is the set of all triples $(a_0,a_1,a_2)$ that are NOT strongly-good and such that
$a_0=0$ if $\nonarrow\not\in\Pset$,
$a_1=0$ if $\{\rightarrow,\leftarrow\}\cap\Pset=\emptyset$,
$a_2=0$ if $\unarrow\not\in\Pset$. The \textindef{strong directed chromatic number}, $\schidir{\Pset}(\hh)$,  of $\hh$ with respect to a palette $\Pset$ is the maximum $\ell +1$ such that for some non-negative integers $a_0,a_1,a_2$, with $a_0+a_1+a_2=\ell$, the triple $(a_0,a_1,a_2)$ is in the clique spectrum of $\hh$. We merely use $\schidir{}(\hh)$ for the strong directed chromatic number if the palette is understood.~\\
\end{defn}

\begin{rem}~\\
\begin{itemize}
   \item The weak [strong] clique spectrum with respect to a given palette is again a downset in the partially ordered set of $r$-tuples ordered coordinatewise.  That is, if $(a_0,a_1,a_2)$ is in the weak [strong] clique spectrum and $(a_0',a_1',a_2')$ has the property that $0\leq a_i'\leq a_i$ for $i=0,1,2$, then $(a_0',a_1',a_2')$ is also in that weak [strong] clique spectrum.
   \item For any palette $\Pset$ and any hereditary property of digraphs, $\hh$, $\wchidir{\Pset}(\hh)\leq\schidir{\Pset}(\hh)$.
   \item If the palette $\Pset\in\{\Pset_{\rm undir},\Pset_{\rm tourn}\}$, the weak and strong directed chromatic numbers are equal and so in those cases, we can use $\chidir{\Pset}=\wchidir{\Pset}=\schidir{\Pset}$.
   \item If the palette is $\Pset_{\rm undir}=\{\nonarrow,\unarrow\}$, then $\chidir{}(\hh)$ is both the binary chromatic number of hereditary property $\hh$.
   \item If $\hh=\forb(H)$ and the palette is $\Pset_{\rm tourn}=\{\leftarrow,\rightarrow\}$, the case of tournaments, then $\chidir{}(H)$ is the fewest number of transitive subtournaments into which $V(H)$ can be partitioned.~\\
%\footnote{Added this remark and removed this: "Note that monotonicity holds for the clique spectrum here as well.  If $(a_0,a_1,a_2)$ is in a clique spectrum and $a_i'\leq a_i$ for $i=0,1,2$, then $(a_0',a_1',a_2')$ is in that clique spectrum."}
\end{itemize}
\end{rem}

%For a hereditary property $\hh=\bigcap_{H\in\F(\hh)}\forb(H)$, and a palette $\Pset$ we define the \textdef{clique spectrum of $\hh$ with respect to $\Pset$}\footnote{define this in $r$-sets? or comment that it can be extended?} is the set of all triples $(a_0,a_1,a_2)$ such that (1) $a_0=0$ if $\nonarrow\not\in\Pset$, $a_1=0$ if $\{\leftarrow,\rightarrow\}\not\subseteq\Pset$, $a_2=0$ if $\unarrow\not\in\Pset$ and (2) for any $H\in\F(\hh)$, then $H$ cannot be partitioned into $\ell$ subsets such that the initial set of $a_0$ of them have no edges labeled $\nonarrow$, the next $a_1$ of them are have no edges labeled $\leftarrow$ or $\rightarrow$ and the final $a_2$ of them are cliques have no edges labeled $\unarrow$.
%
%With that definition, we define the \textdef{$\Pset$-directed chromatic number of $\hh$, $\chidir^{\Pset}(\hh)$}, to be the maximum $\ell+1$ for which there exists a vector $(a_0,a_1,a_2)$, in the clique spectrum of $\hh$ with respect to $\Pset$, such that $\sum_{i=0}^2a_i=\ell$.
%

\subsection{A simple editing algorithm}~\\
\label{sec:digraphs:simplealg}
Let $\Pset$ be a palette and let $\hh$ be a hereditary property of digraphs such that $\hh = \bigcap_{H\in\F(\hh)}\forb(H)$ and each edge of each $H\in\F(\hh)$ has a color in $\Pset$. Further, let $\ell=\wchidir{\Pset}(\hh)-1$ and $(a_0,a_1,a_2)$ be in the weak clique spectrum and $\sum_{i=0}^2a_i=\ell$.  Recall that if a color is not in the palette, then its corresponding $a_i$ value must be set to zero.
% and no digraph in $\F(\hh)$ can be partitioned into $\ell$ subsets such that $a_0$ of them are independent sets, $a_1$ of them are transitive subtournaments and $a_2$ of them are cliques with undirected edges.

%In order for us to bring a given digraph $G$ to a digraph in $\hh$ via edge recolorings, it is sufficient to bring it to a digraph consisting of $\ell$ sets whose edges are specified by the vector $(a_0,a_1,a_2)$. This digraph contains none of the forbidden digraphs.

Partition $V$ into $3$ sets, $S_0,S_1,S_2$ and further refine the partition such that $S_i=V_{i,1}\cup\cdots\cup V_{i,a_i}$, for $i=0,1,2$ and then recolor the edges induced by each $V_{i,j}$ as follows:
\begin{itemize}
   \item If $i=0$, then recolor the edges colored $\nonarrow$ into some other arbitrary color in the palette.
   \item If $i=1$, then recolor the edges $\leftarrow$ and $\rightarrow$ so that there are no directed cycles among those directed edges.
   \item If $i=2$, then recolor the edges colored $\unarrow$ into some other arbitrary color in the palette.
\end{itemize}
This new coloring does not contain any $H\in\F(\hh)$, otherwise the triple $(a_0,a_1,a_2)$ would be weakly good for some $H$. As in the multicolor case, if the partition into sets $V_{i,j}$ is an equipartition, then
$$ \dist(\hh)\leq\frac{1}{\ell}=\frac{1}{\wchidir{\Pset}(\hh)-1} . $$~\\

\subsection{Main results}~\\
In Section~\ref{sec:multicol}, we have seen a general bound in the $r$-graph case. Here, we show a similar result in the directed case.
\begin{theorem}\label{thm:digraphs:easybounds}
Let $\Pset$ be a palette and $\hh$ be a hereditary property of digraphs.
Let $\wchidir{\Pset}=\wchidir{\Pset}(\hh)$  and $\schidir{\Pset}=\schidir{\Pset}(\hh)$ be the weak and strong directed chromatic numbers, respectively, of $\hh$. Recall that if $\Pset\in\{\Pset_{\rm undir},\Pset_{\rm tourn}\}$, then $\chidir{\Pset}=\wchidir{\Pset}=\schidir{\Pset}$. Then,
\begin{enumerate}
   \setcounter{enumi}{-1}
   \item $\ds\frac{1}{4(\schidir{\Pset}-1)}\leq\dist(\hh)\leq\frac{1}{\wchidir{\Pset}-1}$, if $\Pset=\Pset_0=\{\nonarrow,\leftarrow,\rightarrow,\unarrow\}$.
   \item $\ds\frac{1}{3(\schidir{\Pset}-1)}\leq\dist(\hh)\leq\frac{1}{\wchidir{\Pset}-1}$, if $\Pset=\Pset_{\rm compl}=\{\unarrow,\leftarrow,\rightarrow\}$.
   \item $\ds\frac{1}{3(\schidir{\Pset}-1)}\leq\dist(\hh)\leq\frac{1}{\wchidir{\Pset}-1}$, if $\Pset=\Pset_{\rm orien}=\{\nonarrow,\leftarrow,\rightarrow\}$.
   \item $\ds\frac{1}{2(\chidir{\Pset}-1)}\leq\dist(\hh)\leq\frac{1}{\chidir{\Pset}-1}$, if $\Pset=\Pset_{\rm undir}=\{\nonarrow,\unarrow\}$. \label{it:digraphs:easybounds:undir}
   \item $\ds\dist(\hh)=\frac{1}{2(\chidir{\Pset}-1)}$, if $\Pset=\Pset_{\rm tourn}=\{\leftarrow,\rightarrow\}$. \label{it:digraphs:easybounds:tourn}
\end{enumerate}
\end{theorem}

%The bounds in this theorem are tight.\footnote{Once again, I can't prove that.}
We prove Theorem~\ref{thm:digraphs:easybounds} in Section~\ref{sec:digraphs:proofs}. As in the multicolor case, the upper bound is a consequence of the simple editing algorithm.  The lower bound comes from Theorem~\ref{thm:digraphs:basics}, stated below, which is the digraph version of Theorem~\ref{thm:multicol:basics} and deals with computing the edit distance for given hereditary properties of digraphs. In order to do so, we need to investigate the so-called edit distance function, which computes the edit distance of a digraph such that nonedges, directed edges and undirected edges having a specified density.~\\

\subsection{The edit distance function}~\\
\subsubsection{Preliminary definitions}
For a digraph, $G=(V,c)$ with $c:(V)_2\rightarrow\Aset$ and $c$ having the required symmetries as in Definition~\ref{defn:digraphs:basic}, partition $(V)_2$ as follows:
\begin{itemize}
   \item $E_{\nonarrow}(G)$ is the set of all unordered pairs $\{v,w\}$ such that $c(v,w)=\nonarrow$,
   \item $E_{\leftarrow}(G)$ is the set of all ordered pairs $(v,w)$ such that $c(v,w)=\leftarrow$,
   \item $E_{\rightarrow}(G)$ is the set of all ordered pairs $(v,w)$ such that $c(v,w)=\rightarrow$,
   \item $E_{\unarrow}(G)$ is the set of all unordered pairs $\{v,w\}$ such that $c(v,w)=\unarrow$,
\end{itemize}
The definition of a density vector in the $r$-graph case does not translate well to the directed case because of the asymmetry that results from directed edges, so we have a new definition.

Given a palette, $\Pset$, A \textdef{directed density vector $(p,q)$ with respect to $\Pset$} (or, simply, \textdef{density vector} or \textdef{probability vector} where the context is understood) is a nonnegative real vector with the property that $p+2q\leq 1$. Furthermore,
\begin{enumerate}
   \item If $\Pset=\Pset_{\rm compl}=\{\unarrow,\leftarrow,\rightarrow\}$, then $p+2q=1$.
   \item If $\Pset=\Pset_{\rm orien}=\{\nonarrow,\leftarrow,\rightarrow\}$, then $p=0$ and $q\leq 1/2$.
   \item If $\Pset=\Pset_{\rm undir}=\{\nonarrow,\unarrow\}$, then $q=0$ and $p\leq 1$. This is the $r$-graph case where $r=2$ or simply the case of undirected graphs.  See~\cite{AS1} and~\cite{AKM}.
   \item If $\Pset=\Pset_{\rm tourn}=\{\leftarrow,\rightarrow\}$, then $p=0$ and $1-p-2q=0$, so $q=1/2$.
\end{enumerate}
For any density vector $(p,q)$, and an integer $n$, we denote\footnote{Formally, the sizes of the partitions of the edge set should be integral, so we can take the floor function for the sizes of, say $E_{\unarrow}$, $E_{\leftarrow}$, $E_{\rightarrow}$ and the size of $E_{\nonarrow}$ is what remains. Since we fix $p$ and $q$ and let $n$ approach infinity, this will make no appreciable difference.}
$$ \dist_n((p,q),\hh)=\max\left\{\dist(G,\hh) : \begin{array}{l} |V(G)|=n, |E_{\unarrow}(G)|=p{\textstyle\binom{n}{2}}, |E_{\rightarrow}(G)|=q{\textstyle\binom{n}{2}}, \\
|E_{\leftarrow}(G)|=q{\textstyle\binom{n}{2}}\mbox{ and }|E_{\nonarrow}(G)|=(1-p-2q){\textstyle\binom{n}{2}}\end{array}\right\} . $$
Observe that there are, in fact, four densities here; two are equal and all sum to one. Thus, we only need two parameters.  We choose parameter names as above because the case of $q=0$ gives the classical case of undirected graphs, as we see below.
Later in the paper, we show that the following limit exits, which we call the \textdef{edit distance function}:
$$ \dist((p,q), {\hh})=\lim_{n\rightarrow\infty}\dist_n((p,q),\hh). $$

Having the edit distance function, we see that $\dist(\hh)=\max_{(p,q)}\dist((p,q),\hh)$, where the maximum is taken over all density vectors that are valid under the conditions imposed by the palette.

\subsubsection{Types of colorings}
\label{sec:digraphs:types}

In Section~\ref{sec:digraphs:algorithm-analysis}, we define two functions which are described in terms of dir-types, which allow us to compute
the edit distance function. Later in the paper, we shall provide algorithms for such computing.
%We use standard notation $\J$, $\Z$, $\one$, $\zero$ for the the matrix with all entries one, all entries zero, for the vector with all entries one, all entries zero, respectively.
%For a real vector ${\bf v}$, we write that ${\bf v} \geq 0$ if each component of ${\bf v}$ is non-negative.

\begin{defn}
For a palette $\Pset$, a \textindef{$\Pset$-dir-type} (or \textindef{dir-type} or \textindef{type}, where the context and the palette are understood), $K$, is a pair $(U,\phi)$, where $U$ is a finite set of vertices and $\phi : U\times U \rightarrow 2^{\Pset}\setminus \emptyset$, such that
\begin{itemize}
   \item for distinct $x,y$ and $a\in\{\nonarrow,\unarrow\}$, $\phi(x,y)\ni a$ if and only if $\phi(y,x)\ni a$ and
   \item for distinct $x,y$, $\phi(x,y)\ni\rightarrow$ if and only if $\phi(y,x)\ni\leftarrow$ and
   \item $\phi(x,x)\neq\Pset$.~\footnote{Note that it is possible that $\left|\{\leftarrow,\rightarrow\}\cap\phi(x,x)\right|=1$.}
\end{itemize}
The \textindef{sub-dir-type} of $K$ induced by $W\subseteq U$ is the dir-type achieved by deleting the vertices $U-W$ from $K$.

We say that a digraph $H=(V,c)$ \textindef{embeds in type $K=(U,\phi)$} if there is a map $\gamma:V\rightarrow U$ such that for all vertices $v\neq v'$,
\begin{itemize}
   \item if $\gamma(v)\neq\gamma(v')$, then $c(v,v')\in\phi(\gamma(v),\gamma(v'))$,
   \item if $c_0\in\{\nonarrow,\unarrow\}$ and $c_0\not\in\phi(u,u)$, then $\gamma^{-1}(u)$ has no pair with color $c_0$,
   \item if $\{\leftarrow,\rightarrow\}\cap\phi(u,u)=\emptyset$, then $\gamma^{-1}(u)$ has no directed edge, and
   \item if $\left|\{\leftarrow,\rightarrow\}\cap\phi(u,u)\right|=1$, then $\gamma^{-1}(u)$ has no directed cycle.
\end{itemize}
In other words, there is a mapping $\gamma$ that brings each edge of color $c_0$ to a vertex or an edge containing $c_0$ in its color set, except that if a vertex contains exactly one of $\{\leftarrow,\rightarrow\}$ then the pre-image of that vertex can be ordered transitively with respect to the oriented edges. If $H$ embeds in type $K$, we write $H\arrows K$, otherwise we write $H \not\arrows K$.
For every hereditary property $\hh$, we let $\K(\hh)$ be the set of all dir-types such that none of $\F(\hh)$ embeds in that type, i.e.,
$$ \K(\hh)=\left\{K : K\mbox{ is a dir-type}, H\not\arrows K, \forall H\in\F(\hh)\right\} . $$

%A digraph $G'=(V,c)$ \textindef{has type $K=(\{u_1,\ldots,u_k\},\phi)$} if there is a partition of the vertices of $G'$, $V=V_1\cup\cdots\cup V_k$, such that the colors that are absent in vertices and edges of $K$ force the absence of colors in the sets $V_i$ and pairs of sets $(V_i,V_j)$.  Formally,
%\begin{itemize}
%   \item if $i\neq j$, $v_i\in V_i$ and $v_j\in V_j$, then $c(v_i,v_j)\in\phi(u_i,u_j)$,
%   \item if $v_i,v_i'\in V_i$, $c_0\in\{\nonarrow,\unarrow\}$ and $c_0\not\in\phi(u_i,u_i)$, then $c_0\neq c(v_i,v_i')$,
%   \item if $\{\leftarrow,\rightarrow\}\cap\phi(u_i,u_i)=\emptyset$, then $\phi^{-1}(u_i,u_i)$ has no directed edge, and
%   \item if $\left|\{\leftarrow,\rightarrow\}\cap\phi(u_i,u_i)\right|=1$, then $\phi^{-1}(u_i,u_i)$ has no directed cycle.
%\end{itemize}

We say that an digraph $G'=(V,c)$ \textindef{has type $K=(U,\phi)$ if} $G'$ embeds into $K$ with mapping $\gamma: V\rightarrow U$ and $\gamma$ is surjective.
\end{defn}

We have Fact~\ref{fact:digraphs:embed}, also similar to the $r$-graph case, which generalizes the ideas underlying the simple editing algorithm in Section~\ref{sec:digraphs:simplealg}.
\begin{fact}\label{fact:digraphs:embed}
   If $K$ is a dir-type, $G'$ is of type $K$ and $H$ does not embed into $K$, then $H\not\subseteq G'$.
\end{fact}~\\

\subsection{Editing algorithm using types}~\\
\label{sec:digraphs:algorithm}
Let ${\bf w}=(w_1,\ldots,w_k)$ be a density vector and let $(p_{\nonarrow},p_{\leftarrow},p_{\rightarrow},p_{\unarrow})$ be a density vector.  This latter vector will represent a vector of densities.  The number of ordered pairs $(x,y)$ with color ``$\unarrow$'' will be $p_{\unarrow}(n)_2$ and the number of ordered pairs with color ``$\nonarrow$'' will be $p_{\nonarrow}(n)_2$.  The number of ordered pairs with color ``$\leftarrow$'' is $p_{\leftarrow}(n)_2$ and the number of \textit{ordered} pairs with color ``$\rightarrow$'' is $p_{\rightarrow}(n)_2$.  Consequently, $p_{\unarrow}+p_{\nonarrow}+p_{\leftarrow}+p_{\rightarrow}=1$.

The vector ${\bf w}$ will represent a vector of weights, assigned to the vertices of an dir-type with vertices $u_1,\ldots,u_k$, respectively.

Let $\Pset\subseteq\{\nonarrow,\leftarrow,\rightarrow,\unarrow\}$ be a palette, $\hh$ be a hereditary property and $G=(V,c)$ be a digraph in $\Pset$ such that the density vector is $(p_{\nonarrow},p_{\leftarrow},p_{\rightarrow},p_{\unarrow})$. In order to find an upper bound on $\dist(G, \hh)$, it is sufficient to change $G$ to a digraph such that, for all $H\in\F(\hh)$, $H$ does not embed into the new coloring. In particular, if the resulting coloring has type $K \in \K(\hh)$, then this coloring is in $\hh$.

\begin{algorithm}\label{algor:digraphs:recolor}
Fix such a $K=(U, \phi) \in \K(\hh)$ and try to bring $G$ to a coloring of type $K$ by edge-recoloring. Let $U= \{u_1, \ldots, u_k\}$. Partition the vertices of $G$ randomly into sets $V_1, \ldots, V_k$ such that the probability of a vertex to be in a part $V_i$ is $w_i$. With an ordering of the vertices of $G$ and vertices $x<y$, consider an edge $(x,y)$ of $G$, let $x\in V_i$, $y\in V_j$, for $i,j\in\{1,\ldots,k\}$. If $i\neq j$ and $c(x,y)\not\in\phi(u_i,u_j)$, recolor $(x,y)$ with a color from $\phi(u_i,u_j)$.

Next, consider the edges in $V_i$.  If $\phi(u_i,u_i)$ contains exactly one of $\{\leftarrow,\rightarrow\}$, then consider a random order of the vertices of $V_i$, call it $\sigma$.  Let $x<y$ and both in $V_i$.  If $c(x,y)=\leftarrow$, then recolor $(x,y)$ if and only if $\sigma(x)<\sigma(y)$. If $c(x,y)=\rightarrow$, then recolor $(x,y)$ if and only if $\sigma(x)>\sigma(y)$.   Note that this forces $V_i$ to have no directed cycles. If $\phi(u_i,u_i)\not\ni a$ for some $a\in\{\nonarrow,\unarrow\}$, then recolor any edge with color $a$ to a color in $\phi(u_i,u_i)$. This concludes the algorithm.
\end{algorithm}

Algorithm~\ref{algor:digraphs:recolor} is simply a directed graph version of Algorithm~\ref{algor:multicol:recolor}.  We only needed to address the editing of oriented edges.~\\

\subsubsection{Analysis of the editing algorithm}\label{sec:digraphs:algorithm-analysis}
Let us first consider a pair $(x,y)$. If $c(x,y)\in\{\nonarrow,\unarrow\}$, then the probability that the color of $(x,y)$ is unchanged is
$$ \sum_{1\leq i,j\leq k}w_iw_j{\bf 1}_{c(x,y)\in\phi(u_i,u_j)} . $$
If $c(x,y)\in\{\leftarrow,\rightarrow\}$, then the probability that the color of $(x,y)$ is unchanged is
$$ \sum_{1\leq i<j\leq k}w_iw_j{\bf 1}_{\rightarrow\in\phi(u_i,u_j)} +\sum_{i=1}^kw_i^2\frac{|\{\leftarrow,\rightarrow\}\cap\phi(u_i,u_i)|}{2} = \sum_{1\leq i,j\leq k}w_iw_j\frac{1}{2}\left|\{\leftarrow,\rightarrow\}\cap\phi(u_i,u_j)\right| . $$
It doesn't matter whether we consider the pair $(u_i,u_j)$ or $(u_j,u_i)$ in the last term because\linebreak $|\{\leftarrow,\rightarrow\}\cap\phi(u_i,u_j)|$ is invariant whether $i<j$ or $i>j$.

Now, the expected number of changes is
\begin{eqnarray*}
\E[\#\mbox{ changes}] & = & \binom{n}{2}-\sum_{x,y\in V,~x<y} \Pr((x,y)\mbox{ is not changed}) \\
& = & \binom{n}{2}-\sum_{1\leq i,j\leq k}w_iw_jp_{\nonarrow}\binom{n}{2}{\bf 1}_{\nonarrow\in\phi(u_i,u_j)}-\sum_{1\leq i,j\leq k}w_iw_jp_{\unarrow}\binom{n}{2}{\bf 1}_{\unarrow\in\phi(u_i,u_j)} \\ & & -\sum_{1\leq i,j\leq k}w_iw_j\frac{p_{\leftarrow}+p_{\rightarrow}}{2}\binom{n}{2}|\{\leftarrow,\rightarrow\}\cap\phi(u_i,u_j)| \end{eqnarray*}

Let $p=p_{\unarrow}$, $q=\frac{p_{\leftarrow}+p_{\rightarrow}}{2}$ and so $1-p-2q=p_{\nonarrow}$. For $K=(U,c)$, and $\rho\in\{\nonarrow,\unarrow\}$, the matrix ${\bf A}_{\rho}$ is such that the $(i,j)^{\rm th}$ entry is $1$ if $c(u_i,u_j)\ni\rho$ and zero otherwise.  The matrix ${\bf A}_{\rightarrow}$ is a $\{0,1\}$-matrix with the property that
$$ \left({\bf A}_{\rightarrow}\right)_{ij}=\left|\{\leftarrow,\rightarrow\}\cap c(u_i,u_j)\right| . $$
With $\J$ denoting the $k\times k$ all-ones matrix, then we define
$$ \M_K(\p)=\J-(1-p-2q){\bf A}_{\nonarrow}-p{\bf A}_{\unarrow}-q{\bf A}_{\rightarrow} . $$
Consequently, if ${\bf w}=(w_1, \ldots, w_k)$, then
$\E[\#\mbox{ changes}]={\bf w}^T\M_K(\p){\bf w}\binom{n}{2}$.

As in the $r$-graph case, we define two functions in terms of the matrix $\M_K(\p)$:
\begin{itemize}
   \item $f_K(\p)=\left(\frac{1}{k}\one\right)^T\M_K(\p)\left(\frac{1}{k}\one\right)$ and
   \item $g_K(\p)=\min\left\{{\bf w}^T\M_K(\p){\bf w} : {\bf w}^T\one=1, {\bf w}\geq\zero\right\}$.
\end{itemize}

\begin{note}
In the directed case, each ordered pair can receive one of 4 directions, but the density vectors only have two entries rather than three. This is because the above computation shows that an upper bound on editing any digraph is determined not by the pair $(p_{\leftarrow},p_{\rightarrow})$ but only by $q=(p_{\leftarrow}+p_{\rightarrow})/2$. It is straightforward, by the same arguments as in the proof of Theorem~\ref{thm:multicol:basics}, to see that the lower bound for the maximum edit distance is asymptotically achieved by a random graph in which the probability of a forward arc is equal to the probability of a backward arc.
\end{note}~\\

\subsection{Basic results on digraphs}~\\
Theorem~\ref{thm:digraphs:basics} is a parallel to Theorem~\ref{thm:multicol:basics} and summarizes some facts about the edit distance function.  Recall that, depending on the palette, there may be further restrictions on the density vector other than the necessary $p+2q\leq 1$. The dimension, $r$, of the palette, $\Pset$, is the number of members of $\{\nonarrow,\rightarrow,\unarrow\}$ that $\Pset$ has.

\begin{theorem}\label{thm:digraphs:basics}
Let $\hh$ be a hereditary property of digraphs and $\Pset$ a palette. Fix a density vector with respect to $\Pset$, $\p=(p,q)$. The limit $\dist(\p,{\hh}):=\lim_{n\rightarrow\infty}\dist_n(\p,\hh)$ exists. Moreover,
\begin{enumerate}
\item $\dist(\p,\hh) = \inf_{K\in\K(\hh)}f_K(\p) = \inf_{K\in\K(\hh)}g_K(\p)$; \label{it:digraphs:inf}
\item Fix $\epsilon>0$, then with probability approaching 1 as $n\rightarrow\infty$,
$$ \dist(\p, \hh)-\epsilon\leq\dist(G(n,\p),\hh) \leq \dist(\p, \hh) ; $$
\item $\dist(\p, \hh)=\lim_{n\rightarrow\infty}{\rm\bf E}[\dist(G(n,\p),\hh)]$;
\item $\dist(\p, \hh)$ is continuous over the domain of density vectors with respect to $\Pset$ and is concave down;\label{it:digraphs:concon}
\item $\dist(\p,\hh)$ achieves its maximum, $\dist({\hh})$, at some density vector $\p^*_{\hh}$ (in fact, denote the set of all such vectors $\p^*_\hh$) and so,
$$ \dist(\hh)=\lim_{n\rightarrow\infty}{\rm\bf E}[\dist(G(n,\p^*_{\hh}),\hh)] \mbox{; and} $$
\item Both $\p_{\hh}^*$ and $\dist(\hh)$ exist and $\p_{\hh}^*$ is a convex and closed set in $[0,1]^{r-1}$.
\end{enumerate}
\end{theorem}

\begin{note}
Again, we abuse notation so that $\p^*_\hh$ can be a single vector or a set.
\end{note}~\\

\subsection{Example: tournaments}~\\
The case of tournaments is relatively straightforward.  Because in tournaments, there are no edges labeled $\nonarrow$ or $\unarrow$, there is only one density vector, $\p=(0,1/2)$. This means that we only need to consider tournaments that are random, that each arc is forward independently with probability $1/2$.  This leads to a rather simple expression for the edit distance:
\begin{thm}\label{thm:digraphs:tournaments} Let $\hh$ be a nontrivial hereditary property of tournaments and let $\Pset=\Pset_{\rm tourn}=\{\leftarrow,\rightarrow\}$.  Then,
$$ \dist(\hh)=\frac{1}{2(\chidir{\Pset}(\hh)-1)} . $$
\end{thm}

Note that in the case of tournaments, the directed chromatic number of tournament $H$, $\chidir{\Pset_{\rm tourn}}(H)$ is the smallest number of transitive subtournaments into which $H$ can be partitioned. We prove Theorem~\ref{thm:digraphs:tournaments} in Section~\ref{sec:digraphs:tournaments}.~\\

\subsection{Example: triangles}~\\
Theorem~\ref{thm:digraphs:triangles} gives some basic results on examples of hereditary properties of digraphs defined by triangles.  The proof is in Section~\ref{sec:digraphs:triangles}.
\begin{thm}\label{thm:digraphs:triangles}
Consider hereditary properties of digraphs.
\begin{enumerate}
\item If $\F$ is a family that consists of a single directed triangle, then, regardless of the palette, $\dist(\Forb(\F))=1/2$.\label{it:digraphs:triangles:1}
\item If $\F$ is a family that consists of a single transitive triangle and $\Pset=\Pset_{\rm tourn}$, the palette of tournaments, then $\forb(\F)$ is a trivial hereditary property.\label{it:digraphs:triangles:2}
\item If $\F$ is a family of that consists of a single transitive triangle, then, if $\Pset$ is any palette other than $\Pset_{\rm tourn}$, then $\dist(\Forb(\F))=1/2$.\label{it:digraphs:triangles:3}
\item If $\F$ is a family that consists of both a transitive and a directed triangle, and $\Pset$ is any palette other than $\Pset_{\rm tourn}$, then $\dist(\Forb(\F))=1/2$.\label{it:digraphs:triangles:4}
\end{enumerate}
\end{thm}~\\

\subsection{Proofs}~\\
\label{sec:digraphs:proofs}
\subsubsection{Proof of Theorem \ref{thm:digraphs:easybounds}}
\label{sec:digraphs:easyboundproof}
The upper bound for this theorem is proven by the simple editing algorithm from Section~\ref{sec:digraphs:simplealg}.

Let $r=|\Pset|$. For the lower bound, we apply part \ref{it:digraphs:inf} of Theorem~\ref{thm:digraphs:basics}, which states that $\dist(\p,\hh)=\inf_{K\in\K(\hh)}f_K(\p)$. Consider an arbitrary $K=(V,\phi)\in\K(\hh)$, a $\Pset$-dir-type on $k$ vertices. Let $\tilde{K}$ be a graph with vertex set $V$ such that $v$ and $v'$ are adjacent in $\tilde{K}$ if and only if $\phi(v,v')=\Pset$. We observe that $\tilde{K}$ has no clique on $\schidir{\Pset}$ vertices, otherwise for some $H\in\F(\hh)$, $H\mapsto K$. Using Tur\'an's theorem, the number of edges of $\tilde{K}$ is at most $\frac{\schidir{\Pset}-2}{\schidir{\Pset}-1}\cdot\frac{k^2}{2}$.
Let $\p=\frac{1}{r}{\bf 1}$. Consider the matrix ${\bf M}_K(\p)$ and observe that every entry is either zero or is a positive integer multiple of $1/r$. The zero entries correspond exactly to pairs with $\phi$ value equal to $\Pset$. Thus, this matrix ${\bf M}_K(\p)$ has at least $k^2-2\left(\frac{\schidir{\Pset}-2}{\schidir{\Pset}-1}\cdot\frac{k^2}{2}\right)\geq \frac{k^2}{\schidir{\Pset}-1}$ entries with value at least $1/r$.
Therefore, $f_K(\p)=\frac{1}{k^2}{\bf 1}^T{\bf M}_K(\p){\bf 1}$ is at least $1/{r(\schidir{\Pset}-1)}$. Since $K$ was arbitrary, this gives a lower bound for $\dist(\p,\hh)$.~\\

\subsubsection{Proof of Theorem \ref{thm:digraphs:basics}}
\label{sec:digraphs:basics}

The proof of most of this theorem is identical to that of Theorem~\ref{thm:multicol:basics}, which is found in Section~\ref{sec:multicol:basics}. The only significant wrinkle is the upper bound.  That is, if $G$ is a digraph with $p\binom{n}{2}$ edges with color $\unarrow$ and $(1-p-2q)\binom{n}{2}$ edges with color $\nonarrow$, then, with $\p=(p,q)$,
$$ \dist(G,\hh)/{\textstyle \binom{n}{2}}\leq \inf_{K\in\K(\hh)}f_K(\p) . $$

This follows directly from the analysis of the editing algorithm using types from Section~\ref{sec:digraphs:algorithm}.~\\
%This is why, even though each ordered pair can receive one of four distinct colors, our density vectors $(p,q)$ will have at most two independent entries.

\subsubsection{Proof of Theorem~\ref{thm:digraphs:tournaments}}~\\
\label{sec:digraphs:tournaments}
In this case, $\p=(0,1/2)$. Let $\hh$ be a hereditary property of tournaments and $\chidir{}=\chidir{\Pset_{\rm tourn}}(\hh)$. In any type $K$ on $k$ vertices, the vertices have color ``$\rightarrow$'' and the edges either have one direction or both.  By the definition of the directed chromatic number, $H\mapsto K$ if $K$ has a clique of order $\chidir{}$ such that every edge of $K$ has color set $\{\leftarrow,\rightarrow\}$.

Similar to the argument in Section~\ref{sec:digraphs:easyboundproof}, we can use Tur\'an's theorem to find a lower bound for $f_K(\p)$. The bilinear form ${\bf 1}^T{\bf M}_K(\p){\bf 1}$ counts $\frac{1}{2}|V(K)|+\frac{1}{2}|E_{\leftarrow}(K)|+\frac{1}{2}|E_{\rightarrow}(K)|$, where $E_{\rho}(K)$ is the set of \textit{ordered} pairs with color $\rho$.  Since $|E_{\{\leftarrow,\rightarrow\}}(K)|+|E_{\leftarrow}(K)|+|E_{\rightarrow}(K)|=k(k-1)$, Tur\'an's theorem gives that $|E_{\{\leftarrow,\rightarrow\}}(K)|\leq\frac{\chidir{}-2}{\chidir{}-1}k^2$. Consequently,
$$ f_K(\p)=\frac{1}{k^2}{\bf 1}^T{\bf M}_K(\p){\bf 1}=\frac{1}{k^2}\left[\frac{1}{2}k+\frac{1}{2}k(k-1)-\frac{\chidir{}-2}{\chidir{}-1}k^2\right]=\frac{1}{2(\chidir{}-1)} . $$
This concludes the proof of Theorem~\ref{thm:digraphs:tournaments}.~\\

\subsubsection{Proof of Theorem~\ref{thm:digraphs:triangles}}~\\
\label{sec:digraphs:triangles}

\noindent\textbf{(\ref{it:digraphs:triangles:1})} As to the upper bound, linearly order the vertices so that the number of backward edges (i.e., pairs $\{v_i,v_j\}$ such that $i<j$ and $c(v_i,v_j)=\leftarrow$) is minimized. A greedy ordering results in at most half of such edges being present. Reorient such edges so that they become forward edges, hence $\dist(\forb(\F))\leq 1/2$.  Note that this corresponds to a $K$ that consists of a single vertex which has color $\rightarrow$.

For the lower bound, consider an arbitrary $K\in\K(\hh)$ with vertex set $\{u_1,\ldots,u_k\}$ and $\p=(0,1/2)$. This means that ${\bf M}_K(\p)={\bf J}-\frac{1}{2}{\bf A}_{\rightarrow}$.  I.e., $\left({\bf M}_K(\p)\right)_{i,j}=1-\frac{1}{2}|c(u_i,u_j)\cap\{\leftarrow,\rightarrow\}|$.

Here we use an approach due to Sidorenko~\cite{Sid}.  See also~\cite{MT,Martin1,MM,Martin2}.  For the optimal solution, ${\bf w}^*$ to the quadratic program $g_K(\p)=\min\left\{{\bf w}^T\M_K(\p){\bf w} : {\bf w}^T\one=1, {\bf w}\geq\zero\right\}$, the vector $\M_K(\p){\bf w}^*$ is a constant vector, equal to $g_K(\p){\bf 1}$.

Observe that there can be no entry $({\bf M}_K(\p))_{ii}=0$ because that means, for the corresponding vertex $u_i$, $c(u_i,u_i)\supseteq\{\leftarrow,\rightarrow\}$ and a directed triangle maps to such a vertex.
%If there is an entry $({\bf M}_K(\p))_{ii}=1/2$, then there can be no entry $({\bf M}_K(\p))_{ij}=0$, because that would mean that vertex $v_i$ has $|c(v_i,v_i)\cap\{\leftarrow,\rightarrow\}|=1$ and there is a vertex $v_j$ such that $c(v_i,v_j)\supseteq\{\leftarrow,\rightarrow\}$. A directed triangle maps to such a pair of vertices with two vertices of the triangle mapping to $v_i$. So, in that case, each entry of the $i$-th row is at least $1/2$ and so $g_K(p)\geq 1/2$.
Suppose there is some entry $({\bf M}_K(\p))_{ij}=0$. This implies that there are a pair of vertices, $u_i$ and $u_j$ such that $c(u_i,u_j)\supseteq\{\leftarrow,\rightarrow\}$. We observe that $|c(u_i,u_i)\cap\{\leftarrow,\rightarrow\}|=|c(u_j,u_j)\cap\{\leftarrow,\rightarrow\}|=0$, otherwise the directed triangle would map to these two vertices of $K$.  Consequently, $({\bf M}_K(\p))_{ii}=({\bf M}_K(\p))_{jj}=1$. Moreover, for every $\ell\in\{1,\ldots,k\}-\{i,j\}$, we have $({\bf M}_K(\p))_{i\ell}+({\bf M}_K(\p))_{j\ell}\geq 1$. If not, then without loss of generality, we have a triangle $\{u_i,u_j,u_{\ell}\}$ in $K$ such that two edges contain $\{\leftarrow,\rightarrow\}$ and the third contains one of $\{\leftarrow,\rightarrow\}$.  It is easy to see that a directed triangle maps to three such vertices.  But then,
\begin{eqnarray*}
   \lefteqn{\sum_{\ell}({\bf M}_K(\p))_{i\ell}\cdot w_{\ell}+\sum_{\ell}({\bf M}_K(\p))_{j\ell}\cdot w_{\ell}} \\
   & \geq & (1-w_i-w_j)+({\bf M}_K(\p))_{ii}\cdot w_i+({\bf M}_K(\p))_{ij}\cdot w_j+({\bf M}_K(\p))_{ji}\cdot w_i+({\bf M}_K(\p))_{jj}\cdot w_j \\
   & = & (1-w_i-w_j)+w_i+0+0+w_j=1 .
\end{eqnarray*}
Since each sum on the left hand side must be equal to $g_K(\p)$, it must be that $g_K(\p)\geq 1/2$.

Finally, if there is no zero entry in the $i$-th row of ${\bf M}_K(\p)$, then $\sum_{\ell}({\bf M}_K(\p))_{i\ell}w_{\ell}\geq 1/2$. Thus, in all cases, $g_K(\p)\geq 1/2$.~\\

\noindent\textbf{(\ref{it:digraphs:triangles:2})} Here we make the easily verified observation that any tournament with at least $4$ vertices has a transitive subtournament of size $3$. So, the hereditary property consists of no tournaments of size $4$ or more.~\\

\noindent\textbf{(\ref{it:digraphs:triangles:3})} As to the upper bound, equipartition the vertex set arbitrarily and recolor an each edge inside either part to have a color other than one in $\{\leftarrow,\rightarrow\}$. Hence, $\dist(\forb(\F))\leq 1/2$. Note that this corresponds to a $K$ that consists of two vertices colored with some nonempty subset of $\{\nonarrow,\unarrow\}$ and an edge colored $\Pset$.

For the lower bound, simply let $G$ be a transitive tournament.  After editing $G$ to make $G'$, there can be no triangles from $G$ that remain and so Mantel's theorem gives that
$$ \dist(G,G')\geq \frac{1}{\binom{n}{2}}\left(\binom{n}{2}-\left\lfloor\frac{n^2}{4}\right\rfloor\right)=\frac{1}{2}-O\left(\frac{1}{n}\right) . $$

%It should be noted at this stage that, in our definition of the density vector, we always ensure that both backward and forward arcs are equally likely, resulting in the single parameter $q$.\footnote{We may want to dive into this a bit more as to why we can choose just one $q$.}~\\

\noindent\textbf{(\ref{it:digraphs:triangles:4})} Here we can use the trivial fact that if $\hh$ is the hereditary property that forbids both directed and transitive triangles and $\hh'$ is the larger hereditary property in~(\ref{it:digraphs:triangles:3}) which forbids only the transitive triangle, then $\dist(\hh)\geq\dist(\hh')=1/2$.  But the example above of a type $K$ that consists of two vertices with none of $\{\leftarrow,\rightarrow\}$ in its color set is in $\K(\hh)$ in this case as well.  Hence, $f_K(\p)=1/2$ and so $\dist(\hh)=1/2$.~\\

\bibliographystyle{plain}

\begin{thebibliography}{99}
\bibitem{AM} M. Axenovich and R. Martin, A version of Szemer\'edi's regularity lemma for multicolored graphs and directed graphs that is suitable for induced graphs, manuscript.
\end{thebibliography}

\end{document}